\documentclass[12pt,a4paper]{article}
\usepackage{amssymb}

\usepackage{graphicx}
\usepackage{amsmath}
\usepackage{makeidx}
\usepackage{indentfirst}

\newtheorem{result}{Result}[section]

\newcounter{resultnum}[section]
\newtheorem{conclusion}{Conclusion}[section]

\newcounter{conclusionnum}[section]

\newcounter{conditionnum}[section]

\newcounter{conjecturenum}[section]
\newtheorem{example}{Example}[section]

\newcounter{examplenum}[section]

\newcounter{exercisenum}[section]
\newtheorem{lemma}{Lemma}[section]

\newcounter{lemmanum}[section]

\newcounter{notationnum}[section]
\newtheorem{theorem}{Theorem}[section]

\newcounter{theoremnum}[section]
\newtheorem{definition}{Definition}[section]

\newcounter{definitionnum}[section]
\newtheorem{corollary}{Corollary}[section]

\newcounter{corollarynum}[section]
\newtheorem{remark}{Remark}[section]

\newcounter{remarknum}[section]
\newtheorem{proposition}{Proposition}[section]

\newcounter{propositionnum}[section]

\newcounter{acknowledgementnum}[section]

\newcounter{algorithmnum}[section]

\newcounter{axiomnum}[section]

\newcounter{casenum}[section]

\newcounter{claimnum}[section]

\newcounter{summarynum}[section]

\newcounter{problemnum}[section]
\newenvironment{proof}[1][]{\textbf{Proof.} }{}

\begin{document}

\title{Nonholonomic Ricci Flows:\ \ I. Riemann \\
Metrics and Lagrange--Finsler Geometry}
\date{February 21, 2007}
\author{ Sergiu I. Vacaru\thanks{%
sergiu$_{-}$vacaru@yahoo.com, svacaru@fields.utoronto.ca } \\
{\quad} \\
\textsl{The Fields Institute for Research in Mathematical Science} \\
\textsl{222 College Street, Second Floor}\\
\textsl{Toronto, Ontario M5T 3J1, Canada} }
\maketitle

\begin{abstract}
In this paper, it is elaborated the theory the Ricci flows for manifolds
enabled with nonintegrable (nonholonomic) distributions defining nonlinear
connection structures. Such manifolds provide a unified geometric arena for
nonholonomic Riemannian spaces, Lagrange mechanics, Fins\-ler geometry, and
various models of gravity (the Einstein theory and string, or gauge,
generalizations). We follow the method of nonhlonomic frames with associated
nonlinear connection structure and define certain classes of nonholonomic
constraints on Riemann manifolds for which various types of generalized
Finsler geometries can be modelled by Ricci flows. We speculate on possible
applications of the nonholonomic flows in modern geometry, geometric
mechanics and physics.

\vskip0.3cm



\end{abstract}

\newpage


\section{ Introduction}

A series of most remarkable results in mathematics are related to Grisha
Perelman's proof of the Poincare Conjecture \cite{per1, per2,per3} built on
geometrization (Thurston) conjecture \cite{thur01,thur02}, for three
dimensional Riemannian manifolds, and R. Hamilton's Ricci flow theory \cite%
{ham1,ham2}, see reviews and basic references in \cite%
{cao,chen,kleiner,rbook}. Much of the works on Ricci flows has been
performed and validated by experts in the area of geometric analysis and
Riemannian geometry.

Some geometric approaches in modern gravity and string theory are connected
to the method of moving frames and distributions of geometric objects on
(semi) Riemannian manifolds and their generalizations to spaces provided
with nontrivial torsion, nonmetricity and/or nonlinear connection structures %
\cite{cartan1,witten}. The geometry of nonholonomic manifolds\footnote{%
a rigorous definition will be presented few paragraphs below, see also
Definition \ref{defanhm}} and non--Riemannian spaces\footnote{%
as particular cases, we can consider the Riemannian--Finsler and
Lagrange--Hamilton geometry, nonholonomic Lie algebroids, Riemann--Cartan
and metric--affine spaces; in brief, all such spaces will be called
non--Riemannian even some Finsler geometries can be equivalently modelled as
Riemannian nonholonomic manifolds} is largely applied in mechanics and
classical/quantum field theory \cite%
{ma1,ma2,mhh,deleon1,deleon2,marsrat,castro1,castro2,gall,aziz1,aziz2,vncg,vsgg,vcla}%
. Such spaces are characterized by three fundamental geometric objects:\
nonlinear connection (N--connection), linear connection and metric. There is
an important geometrical problem to prove the existence of the ''best
possible'' metric and linear connection adapted to a N--connection
structure. From the point of view of Riemannian geometry, the Thurston
conjecture only asserts the existence of a best possible metric on an
arbitrary closed three dimensional (3D) manifold.

It is a very difficult task to define Ricci flows of mutually compatible
fundamental geometric structures on non--Riemannian manifolds (for instance,
on a Finsler manifold). For such purposes, we can also apply the Hamilton's
approach but correspondingly generalized in order to describe nonholonomic
(constrained) configurations. The first attempts to construct exact
solutions of the Ricci flow equations on nonholonomic Einstein and
Riemann--Cartan (with nontrivial torsion) manifolds, generalizing well known
classes of exact solutions in Einstein and string gravity, were performed in
Refs. \cite{vrf,vv1,vv2} (on extracting holonomic solutions see \cite{crvis}%
).

We take a unified point of view to Riemannian and generalized
Finsler--Lagrange spaces following the geometry of nonholonomic manifolds
and exploit the similarities and emphasize differences between locally
isotropic and anisotropic Ricci flows. In our works, it will be shown when
the remarkable Perelman--Hamilton results hold true for more general
non--Riemannian configurations. It should be noted that this is not only a
straightforward technical extension of the Ricci flow theory to certain
manifolds with additional geometric structures. The problem of constructing
the Finsler--Ricci flow theory contains a number of new conceptual and
fundamental issues on compatibility of geometrical and physical objects and
their optimal configurations.

There are at least three important arguments supporting the investigation of
nonholonomic Ricci flows:\ 1) The Ricci flows of a Riemannian metric may
result in a Finsler like metric if the flows are subjected to certain
nonintegrable constraints and modelled with respect to nonholonomic frames
(we shall prove it in this work).\ 2) Generalized Finsler like metrics
appear naturally as exact solutions in Einstein, string, gauge and
noncommutative gravity, parametrized by generic off--diagonal metrics,
nonholonomic frames and generalized connections (see summaries of results
and methods in \cite{vncg,vsgg}). It is an important physical task to
analyze Ricci flows of such solutions as well of other physically important
solutions (for instance, black holes, solitonic and/pp--waves solutions,
Taub NUT configurations \cite{vrf,vv1,vv2}) resulting in nonholonomic
geometric configurations.\ 3) Finally, the fact that a 3D manifold posses a
"best" Riemannian metric, which implies certain fundamental consequences
(for instance) for our spacetime topology, does not prohibit us to consider
other types of "also not bead" metrics with possible local anisotropy and
nonholonomic gravitational interactions. What are the natural evolution
equations for such configurations and how we can relate them to the topology
of nonholonomic manifolds? We shall address such questions in this (for
regular Lagrange systems) and our further works.

The notion of nonholonomic manifold was introduced independently by G. Vr\v{a%
}nceanu \cite{vr1} and Z. Horak \cite{hor} as a need for geometric
interpretation of nonholonomic mechanical systems (see modern approaches,
criticism and historical remarks in Refs. \cite{bejf,vsgg,grozleit}).%
\footnote{%
For simplicity, we assume throughout this article that all manifolds are
smooth and orientable, even they are provided, or not, with nonholonomic
distributions. In literature, it is used the equivalent term
''anholonomic''. One should be noted here that different types of
anholonomic geometries have been elaborated, for different nonintegrable
(nonholonomic) structures on manifolds, by various schools in geometry,
mechanics and field theory which resulted in certain confusion in
terminology and priorities. In order to avoid ambiguities, we ask the reader
to follow our definitions.} A pair $(M,\mathcal{D})$, where $M$ is a
manifold and $\mathcal{D}$ is a nonintegrable distribution on $M$, is called
a nonholonomic manifold. Three well known classes of nonholonomic manifolds,
when the nonholonomic distribution defines a nonlinear connection
(N--connection) structure, are defined by the Finsler spaces \cite%
{cartan,bejancu,bcs} and their generalizations as Lagrange and Hamilton
spaces \cite{kern,ma1,ma2,mhh} (usually such geometries are modelled on the
tangent bundle $TM).$ More recent examples, related to exact off--diagonal
solutions and nonholonomic frames in Einstein/ string/ gauge/ noncommutative
gravity and nonholonomic Fedosov manifolds \cite{vncg,esv,vsgg} also
emphasize nonholonomic geometric structures.

Let us now sketch the Ricci flow program for nonholonomic manifolds and
Lagrange--Finsler geometries. Different models of "locally anisotropic"
spaces can be elaborated for different types of fundamental geometric
structures (metric, nonlinear and linear connections). In general, such
spaces contain nontrivial torsion and nonmetricity fields. It would be a
very difficult technical task to generalize and elaborate new proofs for all
types of non--Riemannian geometries. Our strategy will be different:\ We
shall formulate the criteria when certain type of Finsler like geometries
can be "extracted" (by imposing the corresponding nonholonomic constraints)
from "well defined" Ricci flows of Riemannian metrics. This is possible
because such geometries can be equivalently described in terms of the Levi
Civita connections or by metric configurations with nontrivial torsion
induced by nonholonomic frames. By nonholonomic transforms of geometric
structures, we shall be able to generate certain classes of nonmetric
geometries and/or generalized torsion configurations.

The aim of this paper (the first one in a series of works) is to formulate
the Ricci flow equations on nonholonomic manifolds and prove the conditions
when such configurations (of Finsler--Lagrange type and in modern gravity)
can be extracted from well defined flows of Riemannian metrics and evolution
of preferred frame structures. Further works will be devoted to explicit
generalizations of G. Perelman results \cite{per1,per2,per3} for
nonholonomic manifolds and spaces provided with almost complex structure
generated by nonlinear connections. We shall also construct new classes of
exact solutions of nonholonomic Ricci flow equations, with noncommutative
and/or Lie algebroid symmetry, defining locally anisotropic flows of black
hole, wormhole and cosmological configurations and developing the results
from Refs. \cite{vrf,vv1,vv2,vncg,vsgg,vcla}.

The works is organized as follow:\ One starts with preliminaries on geometry
of nonholonomic manifolds provided with nonlinear connection (N--connection)
structure in Section 2. We show how nonholonomic configurations can be
naturally defined in modern gravity and the geometry of Riemann--Finsler and
Lagrange spaces in Section 3. Section 4 is devoted to the theory of
anholonomic Ricci flows:\ we analyze the evolution of distinguished
geometric objects and speculate on nonholonomic Ricci flows of symmetric and
nonsymmetric metrics. In Section 5, we prove that the Finsler--Ricci flows
can be extracted from usual Ricci flows by imposing certain classes of
nonholonomic constraints and deformations of connections. We also study
regular Lagrange systems and consider generalized Lagrange--Ricci flows. The
Appendix outlines some necessary results from the local geometry of
N--anholonomic manifolds.

\vskip3pt \textbf{Notation Remarks:\ } We shall use both the coordinate free
and local coordinate formulas which is convenient both to introduce compact
denotations and sketch some proofs. The left up/lower indices will be
considered as labels of geometrical objects, for instance, on a nonholonomic
Riemannian of Finsler space. The boldfaced letters will point that the
objects (spaces) are adapted (provided) to (with) nonlinear connection
structure.

\section{Preliminaries:\ Nonholonomic Manifolds}

We recall some basic facts in the geometry of nonholonomic manifolds
provided with nonlinear connection (N--connection) structure. The reader can
refer to \cite{vsgg,vncg,esv,bejf} for details and proofs (for some
important results we shall sketch the key points for such proofs). On
nonholonomic vectors and (co--) tangent bundles and related
Riemannian--Finsler and Lagrange--Hamilton geometries, we send to Refs. \cite%
{ma1,ma2,mhh,bejancu,bcs}.

\subsection{N--connections}

Consider a $(n+m)$--dimensional manifold $\mathbf{V,}$ with $n\geq 2$ and $%
m\geq 1$ (for a number of physical applications, it is equivalently called
to be a physical and/or geometric space). In a particular case, $\mathbf{V=}%
TM,$ with $n=m$ (i.e. a tangent bundle), or $\mathbf{V=E}=(E,M),$ $\dim M=n,$
is a vector bundle on $M,$ with total space $E.$ In a general case, we can
consider a manifold $\mathbf{V}$ provided with a local fibred structure into
conventional ''horizontal'' and ''vertical'' directions. The local
coordinates on $\mathbf{V}$ are denoted in the form $u=(x,y),$ or $u^{\alpha
}=\left( x^{i},y^{a}\right) ,$ where the ''horizontal'' indices run the
values $i,j,k,\ldots =1,2,\ldots ,n$ and the ''vertical'' indices run the
values $a,b,c,\ldots =n+1,n+2,\ldots ,n+m.$\footnote{%
For the tangent bundle $TM,$ we can consider that both type of indices run
the same values.} We denote by $\pi ^{\top }:T\mathbf{V}\rightarrow TM$ the
differential of a map $\pi :\mathbf{V}\rightarrow V$ defined by fiber
preserving morphisms of the tangent bundles $T\mathbf{V}$ and $TM.$ The
kernel of $\pi ^{\top }$ is just the vertical subspace $v\mathbf{V}$ with a
related inclusion mapping $i:v\mathbf{V}\rightarrow T\mathbf{V}.$

\begin{definition}
A nonlinear connection (N--connection) $\mathbf{N}$ on a manifold $\mathbf{V}
$ is defined by the splitting on the left of an exact sequence
\begin{equation*}
0\rightarrow v\mathbf{V}\overset{i}{\rightarrow} T\mathbf{V}\rightarrow T%
\mathbf{V}/v\mathbf{V}\rightarrow 0,
\end{equation*}%
i. e. by a morphism of submanifolds $\mathbf{N:\ \ }T\mathbf{V}\rightarrow v%
\mathbf{V}$ such that $\mathbf{N\circ i}$ is the unity in $v\mathbf{V}.$
\end{definition}

Locally, a N--connection is defined by its coefficients $N_{i}^{a}(u),$%
\begin{equation}
\mathbf{N}=N_{i}^{a}(u)dx^{i}\otimes \frac{\partial }{\partial y^{a}}.
\label{coeffnc}
\end{equation}%
Globalizing the local splitting, one prove:

\begin{proposition}
Any N--connection is defined by a Whitney sum of conventional horizontal (h)
subspace, $\left( h\mathbf{V}\right) ,$ and vertical (v) subspace, $\left( v%
\mathbf{V}\right) ,$
\begin{equation}
T\mathbf{V}=h\mathbf{V}\oplus v\mathbf{V}.  \label{whitney}
\end{equation}
\end{proposition}

The sum (\ref{whitney}) states on $T\mathbf{V}$ a nonholonomic
(equivalently, anholonomic, or nonintegrable) distribution of horizontal and
vertical subspaces. The well known class of linear connections consists on a
particular subclass with the coefficients being linear on $y^{a},$ i.e. $%
N_{i}^{a}(u)=\Gamma _{bj}^{a}(x)y^{b}.$

The geometric objects on $\mathbf{V}$ can be defined in a form adapted to a
N--connection structure, following certain decompositions being invariant
under parallel transports preserving the splitting (\ref{whitney}). In this
case, we call them to be distinguished (by the N--connection structure),
i.e. d--objects. For instance, a vector field $\mathbf{X}\in T\mathbf{V}$ \
is expressed
\begin{equation*}
\mathbf{X}=(hX,\ vX),\mbox{ \ or \ }\mathbf{X}=X^{\alpha }\mathbf{e}_{\alpha
}=X^{i}\mathbf{e}_{i}+X^{a}e_{a},
\end{equation*}%
where $hX=X^{i}\mathbf{e}_{i}$ and $vX=X^{a}e_{a}$ state, respectively, the
adapted to the N--connection structure horizontal (h) and vertical (v)
components of the vector. In brief, $\mathbf{X}$ is called a distinguished
vectors, in brief, d--vector). In a similar fashion, the geometric objects
on $\mathbf{V}$ like tensors, spinors, connections, ... are called
respectively d--tensors, d--spinors, d--connections if they are adapted to
the N--connection splitting (\ref{whitney}).

\begin{definition}
The N--connection curvature is defined as the Neijenhuis tensor,%
\begin{equation}
\mathbf{\Omega }(\mathbf{X,Y})\doteqdot \lbrack vX,vY]+\ v[\mathbf{X,Y}]-v[vX%
\mathbf{,Y}]-v[\mathbf{X,}vY].  \label{njht}
\end{equation}
\end{definition}

In local form, we have for (\ref{njht})
\begin{equation*}
\mathbf{\Omega }=\frac{1}{2}\Omega _{ij}^{a}\ d^{i}\wedge d^{j}\otimes
\partial _{a},
\end{equation*}%
with coefficients%
\begin{equation}
\Omega _{ij}^{a}=\frac{\partial N_{i}^{a}}{\partial x^{j}}-\frac{\partial
N_{j}^{a}}{\partial x^{i}}+N_{i}^{b}\frac{\partial N_{j}^{a}}{\partial y^{b}}%
-N_{j}^{b}\frac{\partial N_{i}^{a}}{\partial y^{b}}.  \label{ncurv}
\end{equation}

Any N--connection $\mathbf{N}$ may be characterized by an associated frame
(vielbein) structure $\mathbf{e}_{\nu }=(\mathbf{e}_{i},e_{a}),$ where
\begin{equation}
\mathbf{e}_{i}=\frac{\partial }{\partial x^{i}}-N_{i}^{a}(u)\frac{\partial }{%
\partial y^{a}}\mbox{ and
}e_{a}=\frac{\partial }{\partial y^{a}},  \label{dder}
\end{equation}%
and the dual frame (coframe) structure $\mathbf{e}^{\mu }=(e^{i},\mathbf{e}%
^{a}),$ where
\begin{equation}
e^{i}=dx^{i}\mbox{ and }\mathbf{e}^{a}=dy^{a}+N_{i}^{a}(u)dx^{i}.
\label{ddif}
\end{equation}%
These vielbeins are called respectively N--adapted frames and coframes. In
order to preserve a relation with the previous denotations \cite{vncg,vsgg},
we emphasize that $\mathbf{e}_{\nu }=(\mathbf{e}_{i},e_{a})$ and $\mathbf{e}%
^{\mu }=(e^{i},\mathbf{e}^{a})$ are correspondingly the former
''N--elongated'' partial derivatives $\delta _{\nu }=\delta /\partial u^{\nu
}=(\delta _{i},\partial _{a})$ and N--elongated differentials $\delta ^{\mu
}=\delta u^{\mu }=(d^{i},\delta ^{a}).$ This emphasizes that the operators (%
\ref{dder}) and (\ref{ddif}) define certain ``N--elongated'' partial
derivatives and differentials which are more convenient for tensor and
integral calculations on such nonholonomic manifolds.\footnote{%
We shall use always ''boldface'' symbols if it would be necessary to
emphasize that certain spaces and/or geometrical objects are
provided/adapted to a\ N--connection structure, or with the coefficients
computed with respect to N--adapted frames.} The vielbeins (\ref{ddif})
satisfy the nonholonomy relations
\begin{equation}
\lbrack \mathbf{e}_{\alpha },\mathbf{e}_{\beta }]=\mathbf{e}_{\alpha }%
\mathbf{e}_{\beta }-\mathbf{e}_{\beta }\mathbf{e}_{\alpha }=W_{\alpha \beta
}^{\gamma }\mathbf{e}_{\gamma }  \label{anhrel}
\end{equation}%
with (antisymmetric) nontrivial anholonomy coefficients $W_{ia}^{b}=\partial
_{a}N_{i}^{b}$ and $W_{ji}^{a}=\Omega _{ij}^{a}.$ The above presented
formulas present the proof of

\begin{proposition}
\label{pddv}A N--connection on $\mathbf{V}$ defines a preferred nonholonomic
N--adapted frame (vielbein) structure $\mathbf{e}=(he,ve)$ and its dual $%
\widetilde{\mathbf{e}}=\left( h\widetilde{e},v\widetilde{e}\right) $ with $%
\mathbf{e}$ and $\widetilde{\mathbf{e}}$ linearly depending on N--connection
coefficients.\
\end{proposition}

For simplicity, we shall work with a particular class of nonholonomic
manifolds:

\begin{definition}
\label{defanhm} A manifold $\mathbf{V}$ is N--anholonomic if its tangent
space $T\mathbf{V}$ it enabled with a N--connection structure (\ref{whitney}%
).
\end{definition}

There are two important examples of N--anholonomic manifolds, when $V=E,$ or
$TM $:

\begin{example}
\label{enhvtb}A vector bundle $\mathbf{E}=(E,\pi ,M,\mathbf{N}),$ defined by
a surjective projection $\pi :E\rightarrow M,$ with $M$ being the base
manifold, $\dim M=n,$ and $E$ being the total space, $\dim E=n+m,$ and
provided with a N--connection splitting (\ref{whitney}) is called
N--anholonomic vector bundle. A particular case is that of N--anholonomic
tangent bundle $\mathbf{TM}=(TM,\pi ,M,\mathbf{N}),$ with dimensions $n=m.$
\end{example}

In a similar manner, we can consider different types of (super) spaces,
Riemann or Riemann--Cartan manifolds, noncommutative bundles, or
superbundles, provided with nonholonomc distributions (\ref{whitney}) and
preferred systems of reference \cite{vncg,vsgg}.

\subsection{Torsions and curvatures of d--connections and d--metrics}

One can be defined N--adapted linear connection and metric structures:

\begin{definition}
\label{ddc}A distinguished connection (d--connection) $\mathbf{D}$ on a
N--anho\-lo\-no\-mic manifold $\mathbf{V}$ is a linear connection conserving
under parallelism the Whitney sum (\ref{whitney}).
\end{definition}

For any d--vector $\mathbf{X,}$ there is a decomposition of $\mathbf{D}$
into h-- and v--covariant derivatives,%
\begin{equation}
\mathbf{D}_{\mathbf{X}}\mathbf{\doteqdot X}\rfloor \mathbf{D=}\ hX\rfloor
\mathbf{D+}\ vX\rfloor \mathbf{D=}Dh_{X}+D_{vX}=hD_{X}+vD_{X}.
\label{dconcov}
\end{equation}
The symbol ''$\rfloor "$ in (\ref{dconcov}) denotes the interior product. We
shall write conventionally that $\mathbf{D=}(hD,\ vD),$ or $\mathbf{D}%
_{\alpha }=(D_{i},D_{a}).$ For convenience, in Appendix, we present some
local formulas for d--connections $\mathbf{D=\{\Gamma }_{\ \alpha \beta
}^{\gamma }=\left( L_{jk}^{i},L_{bk}^{a},C_{jc}^{i},C_{bc}^{a}\right) \},$
with $hD=(L_{jk}^{i},L_{bk}^{a})$ and $vD=(C_{jc}^{i},C_{bc}^{a}),$ see (\ref%
{dcon1}).

\begin{definition}
The torsion of a d--connection $\mathbf{D=}(hD,\ vD)\mathbf{,}$ for any
d--vectors $\mathbf{X,Y}$ is defined by d--tensor field
\begin{equation}
\mathbf{T(X,Y)\doteqdot \mathbf{D}_{\mathbf{X}}Y-D}_{\mathbf{Y}}\mathbf{%
X-[X,Y].}  \label{tors1}
\end{equation}
\end{definition}

One has a N--adapted decomposition
\begin{equation}
\mathbf{T(X,Y)=T(}hX,hY\mathbf{)+T(}hX,\ vY\mathbf{)+T(}vX,hY\mathbf{)+T(}%
vX,\ vY\mathbf{).}  \label{tors2}
\end{equation}

Considering h- and v--projections of (\ref{tors2}) and taking in the account
that $h[vX, vY] =0,$ one proves

\begin{theorem}
\label{tht}The torsion $\mathbf{T}$ of a d--connection $\mathbf{D}$ is
defined by five nontrivial d--tensor fields adapted to the h-- and
v--splitting by the N--connection structure%
\begin{eqnarray*}
hT(hX,hY) &\mathbf{\doteqdot }&D_{hX}\ hY-D_{hY}\ hX-h[\mathbf{X,Y}], \\
vT(hX,hY)&\mathbf{\doteqdot \ }& v[hY,hX], \\
vT(hX, vY) &\mathbf{\doteqdot }&\mathbf{-\ }vD_{vY}\ hX-h[hX, vY], \\
vT(hX,vY) &\mathbf{\doteqdot \ }& vD_{hX}\ vY- v[hX,\mathbf{\ }vY], \\
vT(vX,vY) &\mathbf{\doteqdot }& vD_{X}\ vY - vD_{Y}\ vX-v[vX, vY].
\end{eqnarray*}%
\
\end{theorem}

The d--torsions $hT(hX,hY),vT(vX,vY),...$ are called respectively the $h$ $%
(hh)$--torsion, $v(vv)$--torsion and so on. The local formulas (\ref{dtors})
for torsion $\mathbf{T}$ are given in Appendix.

\begin{definition}
The curvature of a d--connection $\mathbf{D}$ is defined
\begin{equation}
\mathbf{R(X,Y)\doteqdot \mathbf{D}_{\mathbf{X}}\mathbf{D}_{\mathbf{Y}}-D}_{%
\mathbf{Y}}\mathbf{D}_{\mathbf{X}}\mathbf{-D}_{\mathbf{[X,Y]}}  \label{curv1}
\end{equation}%
for any d--vectors $\mathbf{X,Y.}$
\end{definition}

By straightforward calculations, one check the properties%
\begin{eqnarray*}
hR(\mathbf{X,Y})\mathbf{\ }vZ &=&0,\mathbf{\ }vR(\mathbf{X,Y})hZ=0, \\
\mathbf{R(X,Y)Z} &\mathbf{=}&hR\mathbf{(X,Y)}hZ+\mathbf{\ }vR\mathbf{(X,Y)\ }%
vZ,
\end{eqnarray*}%
for any for any d--vectors $\mathbf{X,Y,Z.}$

\begin{theorem}
\label{thr}The curvature $\mathbf{R}$ of a d--connection $\mathbf{D}$ is
completely defined by six d--curvatures
\begin{eqnarray*}
\mathbf{R(}hX\mathbf{,}hY\mathbf{)}hZ &\mathbf{=}&\left(
D_{hX}D_{hY}-D_{hY}D_{hX}-D_{[hX,hY]}-\mathbf{\ }vD_{[hX,hY]}\right) hZ, \\
\mathbf{R(}hX\mathbf{,}hY\mathbf{)\ }vZ &\mathbf{=}&\left(
D_{hX}D_{hY}-D_{hY}D_{hX}-D_{[hX,hY]}-\mathbf{\ }vD_{[hX,hY]}\right) \ vZ, \\
\mathbf{R(}vX\mathbf{,}hY\mathbf{)}hZ &\mathbf{=}&\left( D_{hX}D_{hY}-D_{hY}%
\mathbf{\ }D_{vX}-D_{[\mathbf{v}X,hY]}-vD_{[vX,hY]}\right) hZ, \\
\mathbf{R(}vX\mathbf{,}vY\mathbf{)\mathbf{\ }}vZ &\mathbf{=}&\left( D_{vX}%
\mathbf{\ }D_{hY}-\mathbf{\ }D_{hY}\mathbf{\ }D_{vX}-D_{h[vX,hY]}-D_{v[%
\mathbf{v}X,hY]}\right) vZ, \\
\mathbf{R(}vX\mathbf{,}vY\mathbf{)}hZ &\mathbf{=}&\left(
D_{vX}D_{vY}-D_{vY}D_{vX}-D_{v[vX,\mathbf{\ }vY]}\right) hZ, \\
\mathbf{R(}vX\mathbf{,}vY\mathbf{)}vZ &\mathbf{=}&\left( D_{vX}D_{vY}-D_{vY}%
\mathbf{\ }D_{vX}-D_{v[vX,vY]}\right) vZ.
\end{eqnarray*}
\end{theorem}

The formulas for local coefficients of d--curvatures $\mathbf{R=\{\mathbf{R}%
_{\ \beta \gamma \delta }^{\alpha }\}}$ are given in Appendix, see (\ref%
{dcurv}).

\begin{definition}
A metric structure $\ \breve{g}$ on a N--anholonomic manifold $\mathbf{V}$
is a symmetric covariant second rank tensor field which is not degenerated
and of constant signature in any point $\mathbf{u\in V.}$
\end{definition}

In general, a metric structure is not adapted to a N--connection structure.

\begin{definition}
\label{ddm}A d--metric $\mathbf{g}=hg\oplus _{N}vg$ is a usual metric tensor
which contracted to a d--vector results in a dual d--vector, d--covector
(the duality being defined by the inverse of this metric tensor).
\end{definition}

The relation between arbitrary metric structures and d--metrics is
established by

\begin{theorem}
\label{tdm}Any metric $\ \breve{g}$ can be equivalently transformed into a
d--metric
\begin{equation}
\mathbf{g}=hg(hX,hY)+vg(vX,vY)  \label{dmetra}
\end{equation}%
adapted to a given N--connection structure.
\end{theorem}

\begin{proof}
We introduce denotations $h\breve{g}(hX,hY)=hg(hX,hY)$ and $v\breve{g}(vX,$ $%
vY)=vg(vX,vY)$ and try to find a N--connection when
\begin{equation}
\breve{g}(hX,vY)=0  \label{algn01}
\end{equation}%
for any d--vectors $\mathbf{X,Y.}$ In local form, \ the equation (\ref%
{algn01}) is an algebraic equation for the N--connection coefficients $%
N_{i}^{a},$ see formulas (\ref{metr}) and (\ref{ansatz}) in Appendix. $%
\square $
\end{proof}

\vskip3pt

A distinguished metric (in brief, d--metric) on a N--anholo\-nom\-ic
manifold $\mathbf{V}$ is a usual second rank metric tensor $\mathbf{g}$
which with respect to a N--adapted basis (\ref{ddif}) can be written in the
form%
\begin{equation}
\mathbf{g}=\ g_{ij}(x,y)\ e^{i}\otimes e^{j}+\ h_{ab}(x,y)\ \mathbf{e}%
^{a}\otimes \mathbf{e}^{b}  \label{m1}
\end{equation}%
defining a N--adapted decomposition $\mathbf{g=}hg\mathbf{\oplus _{N}}%
vg=[hg,vg].$

From the class of arbitrary d--connections $\mathbf{D}$ on $\mathbf{V,}$ one
distinguishes those which are metric compatible (metrical d--connections)
satisfying the condition%
\begin{equation}
\mathbf{Dg=0}  \label{metcomp}
\end{equation}%
including all h- and v-projections
\begin{equation*}
D_{j}g_{kl}=0,D_{a}g_{kl}=0,D_{j}h_{ab}=0,D_{a}h_{bc}=0.
\end{equation*}%
Different approaches to Finsler--Lagrange geometry modelled on $\mathbf{TM}$
(or on the dual tangent bundle $\mathbf{T}^{\ast }\mathbf{M,}$ in the case
of Cartan--Hamilton geometry) were elaborated for different d--metric
structures which are metric compatible \cite{cartan,ma1,ma2,mhh} or not
metric compatible \cite{bcs}.

\subsection{(Non) adapted linear connections}

For any metric structure $\mathbf{g}$ on a manifold $\mathbf{V,}$ there is
the unique metric compatible and torsionless Levi Civita connection $\nabla $
for which $\ ^{\nabla }\mathcal{T}^{\alpha }=0$ and $\nabla \mathbf{g=0.}$
This is not a d--connection because it does not preserve under parallelism
the N--connection splitting (\ref{whitney}) (it is not adapted to the
N--connection structure).

\begin{theorem}
\label{thcdc}For any d--metric $\mathbf{g}=[hg,vg]$ on a N--anholonomic
manifold $\mathbf{V,}$ there is a unique metric canonical d--connection $%
\widehat{\mathbf{D}}$ satisfying the conditions $\widehat{\mathbf{D}}\mathbf{%
g=}0$ and with vanishing $h(hh)$--torsion, $v(vv)$--torsion, i. e. $h%
\widehat{T}(hX,hY)=0$ and $\mathbf{\ }v\widehat{T}(vX,\mathbf{\ }vY)=0.$
\end{theorem}

\begin{proof}
By straightforward calculations, we can verify that the d--connec\-ti\-on
with coefficients $\widehat{\mathbf{\Gamma }}_{\ \alpha \beta }^{\gamma
}=\left( \widehat{L}_{jk}^{i},\widehat{L}_{bk}^{a},\widehat{C}_{jc}^{i},%
\widehat{C}_{bc}^{a}\right) ,$ see (\ref{candcon}) in Appendix, satisfies
the condition of Theorem.$\square $
\end{proof}

\begin{definition}
\label{drcm}A N--anholonomic Riemann--Cartan manifold $\ ^{RC}\mathbf{V}$ is
defined by a d--metric $\mathbf{g}$ and a metric d--connection $\mathbf{D}$
structures. For a particular case, we can consider that a space$\ ^{R}%
\widehat{\mathbf{V}}$ is a N--anholonomic Riemann manifold if its
d--connection structure is canonical, i.e. $\mathbf{D=}\widehat{\mathbf{D}}.$
\end{definition}

The d--metric structure $\mathbf{g}$ on$\ ^{RC}\mathbf{V}$ is of type (\ref%
{m1}) and satisfies the metricity conditions (\ref{metcomp}). With respect
to a local coordinate basis, the metric $\mathbf{g}$ is parametrized by a
generic off--diagonal metric ansatz (\ref{ansatz}). For a particular case,
we can take $\mathbf{D=}\widehat{\mathbf{D}}$ and treat the torsion $%
\widehat{\mathbf{T}}$ as a nonholonomic frame effect induced by a
nonintegrable N--splitting. We conclude that a N--anholonomic Riemann
manifold is with nontrivial torsion structure (\ref{dtors}) (defined by the
coefficients of N--connection (\ref{coeffnc}), and d--metric (\ref{m1}) and
canonical d--connection (\ref{candcon})). Nevertheless, such manifolds can
be described alternatively, equivalently, as a usual (holonomic) Riemann
manifold \ with the usual Levi Civita for the metric (\ref{metr}) with
coefficients (\ref{ansatz}). We do not distinguish the existing nonholonomic
structure for such geometric constructions.

For more general applications, we have to consider additional torsion
components, for instance, by the so--called $H$--field in string gravity %
\cite{string1}.

\begin{theorem}
\label{teq}The geometry of a (semi) Riemannian manifold $V$ with prescribed $%
(n+m)$--splitting (nonholonomic h- and v--decomposition) is equivalent to
the geometry of a canonical $^{R}\widehat{\mathbf{V}}.$
\end{theorem}

\begin{proof}
Let $g_{\alpha \beta }$ be the metric coefficients, with respect to a local
coordinate frame, on $V.$ The $(n+m)$--splitting states for a
paramterization of type (\ref{ansatz}) which allows us to define the
N--connection coefficients $N_{i}^{a}$ by solving the algebraic equations (%
\ref{aux1a}) (roughly speaking, the N--connection coefficients are defined
by the ''off--diagonal'' N--coefficients, considered with respect to those
from the blocks $n\times n$ and $m\times m).$ Having defined $\mathbf{N}%
=\{N_{i}^{a}\},$ we can compute the N--adapted frames $\mathbf{e}_{\alpha }$
(\ref{dder}) and $\mathbf{e}^{\beta }$ (\ref{ddif}) by using frame
transforms (\ref{vt1}) \ and (\ref{vt2}) for any fixed values $e_{i}^{\
\underline{i}}(u)$ and $e_{a}^{\ \underline{a}}(u);$ for instance, for
coordinate frames $e_{i}^{\ \underline{i}}=\delta _{i}^{\ \underline{i}}$
and $e_{a}^{\ \underline{a}}=\delta _{a}^{\ \underline{a}}.$ As a result,
the metric structure is transformed into a d--metric of type (\ref{m1}). We
can say that $V$ is equivalently re--defined as a N--anholonomic manifold $%
\mathbf{V.}$

It is also possible to compute the coefficients of canonical d--connection $%
\widehat{\mathbf{D}}$ following formulas (\ref{candcon}). We conclude that
the geometry of a (semi) Riemannian manifold $V$ with prescribed $(n+m)$%
--splitting can be described equivalently by geometric objects on a
canonical N--anholonomic manifold $^{R}\widehat{\mathbf{V}}$ with induced
torsion $\widehat{\mathbf{T}}$ with the coefficients computed by introducing
(\ref{candcon}) into (\ref{dtors}).

The inverse construction also holds true: A d--metric (\ref{m1}) on $^{R}%
\widehat{\mathbf{V}}$ is also a metric on $V$ but with respect to certain
N--elongated basis (\ref{ddif}). It can be also rewritten with respect to a
coordinate bases having the parametrization (\ref{ansatz}). $\square $
\end{proof}

\vskip3pt

From this Theorem, by straightforward computations with respect to
N--adapted bases (\ref{ddif}) and (\ref{dder}), one follows

\begin{corollary}
\label{cteq}The metric of a (semi) Riemannian manifold provided with a
preferred N--adapted frame structure defines canonically two equivalent
linear connection structures: the Levi Civita connection and the canonical
d--connection.
\end{corollary}

\begin{proof}
On a manifold $^{R}\widehat{\mathbf{V}},$ we can work with two equivalent
linear connections. If we follow only the methods of Riemannian geometry, we
have to chose the Levi Civita connection. In some cases, it may be optimal
to elaborate a N--adapted tensor and differential calculus for nonholnomic
structures, i.e. to chose the canonical d--connection. With respect to
N--adapted frames, the coefficients of one connection can be expressed via
coefficients of the second one, see formulas (\ref{lccon}) and (\ref{candcon}%
). Both such linear connections are defined by the same off--diagonal metric
structure. For diagonal metrics with respect to local coordinate frames, the
constructions are trivial.$\square $
\end{proof}

\vskip3pt

Having prescribed a nonholonomic $n+m$ splitting on a manifold $V,$ we can
define two canonical linear connections $\nabla $ and $\widehat{\mathbf{D}}.$
Correspondingly, these connections are characterized by two curvature
tensors, $_{\shortmid }R_{~\beta \gamma \delta }^{\alpha }(\nabla )$
(computed by introducing $_{\shortmid }\Gamma _{\beta \gamma }^{\alpha }$
into (\ref{dconf}) and (\ref{curv})) and $\mathbf{R}_{~\beta \gamma \delta
}^{\alpha }(\widehat{\mathbf{D}})$ (with the N--adapted coefficients
computed following formulas (\ref{dcurv})). Contracting indices, we can
commute the Ricci tensor $Ric(\nabla )$ and the Ricci d--tensor $\mathbf{Ric}%
(\widehat{\mathbf{D}})$ following formulas (\ref{dricci}), correspondingly
written for $\nabla $ and $\widehat{\mathbf{D}}.$ Finally, using the inverse
d--tensor $\mathbf{g}^{\alpha \beta }$ for both cases, we compute the
corresponding scalar curvatures $\ ^{s}R(\nabla )$ and $\ ^{s}\mathbf{R(%
\widehat{\mathbf{D}}),}$ see formulas (\ref{sdccurv}) by contracting,
respectively, with the Ricci tensor and Ricci d--tensor.

\subsection{Metrization procedure and preferred linear connections}

On a N--anholonomic manifold $\mathbf{V},$ with prescribed fundamental
geometric structures $\mathbf{g}$ and $\mathbf{N,}$ we can consider various
classes of d--connections $\mathbf{D,}$ which, in general, are not metric
compatible, i.e. $\mathbf{Dg\neq 0.}$ The canonical d--connection $\widehat{%
\mathbf{D}}$ is the ''simplest'' metric one, with respect to which other
classes of d--connections $\mathbf{D=}$ $\widehat{\mathbf{D}}+\mathbf{Z}$
can be distinguished by their deformation (equivalently, distorsion, or
deflection) d--tensors $\mathbf{Z.}$ Every geometric construction performed
for a d--connection $\mathbf{D}$ can be redefined for $\widehat{\mathbf{D}},$
and inversely, if $\mathbf{Z}$ is well defined.

Let us consider the set of all possible nonmetric and metric d--connections
constructed only form the coefficients of a d--metric and N--connecti\-on
structure, $g_{ij},h_{ab}$ and $N_{i}^{a},$ and their partial derivatives.
Such d--connecti\-ons can be generated by two procedures of deformation,
\begin{eqnarray*}
\widehat{\mathbf{\Gamma }}_{\ \alpha \beta }^{\gamma } &\rightarrow &\ ^{[K]}%
\mathbf{\Gamma }_{\ \alpha \beta }^{\gamma }=\mathbf{\Gamma }_{\ \alpha
\beta }^{\gamma }+\ ^{[K]}\mathbf{Z}_{\ \alpha \beta }^{\gamma }%
\mbox{\
(Kawaguchi's metrization \cite{kaw1,kaw2}) }, \\
\mbox{ or } &\rightarrow &^{[M]}\mathbf{\Gamma }_{\ \alpha \beta }^{\gamma }=%
\widehat{\mathbf{\Gamma }}_{\ \alpha \beta }^{\gamma }+\ ^{[M]}\mathbf{Z}_{\
\alpha \beta }^{\gamma }\mbox{\ (Miron's
connections \cite{ma2} )},
\end{eqnarray*}%
where $\ ^{[K]}\mathbf{Z}_{\ \alpha \beta }^{\gamma }$ and $\ ^{[M]}\mathbf{Z%
}_{\ \alpha \beta }^{\gamma }$ are deformation \ d--tensors.

\begin{theorem}
\ \label{kmp}For given d--metric $\mathbf{g}_{\alpha \beta }=[g_{ij},h_{ab}]$
and N--connection $\mathbf{N}=\{N_{i}^{a}\}$ structures, the deformation
d--tensor \
\begin{eqnarray*}
^{\lbrack K]}\mathbf{Z}_{\ \alpha \beta }^{\gamma } &=&\{\ ^{[K]}Z_{\
jk}^{i}=\frac{1}{2}g^{im}D_{j}g_{mk},\ ^{[K]}Z_{\ bk}^{a}=\frac{1}{2}%
h^{ac}D_{k}h_{cb}, \\
\ \ ^{[K]}Z_{\ ja}^{i} &=&\frac{1}{2}g^{im}D_{a}g_{mj},\ ^{[K]}Z_{\ bc}^{a}=%
\frac{1}{2}h^{ad}D_{c}h_{db}\}
\end{eqnarray*}%
transforms a d--connection $\mathbf{\Gamma }_{\ \alpha \beta }^{\gamma
}=\left( L_{jk}^{i},L_{bk}^{a},C_{jc}^{i},C_{bc}^{a}\right) $\ into a metric
d--connec\-ti\-on%
\begin{equation*}
\ ^{[K]}\mathbf{\Gamma }_{\ \alpha \beta }^{\gamma }=\left( L_{jk}^{i}+\
^{[K]}Z_{\ jk}^{i},L_{bk}^{a}+\ ^{[K]}Z_{\ bk}^{a},C_{jc}^{i}+\ ^{[K]}Z_{\
ja}^{i},C_{bc}^{a}+\ ^{[K]}Z_{\ bc}^{a}\right) .
\end{equation*}%
\
\end{theorem}

\begin{proof}
It consists from a straightforward verification that the conditions
metricity conditions $\ ^{[K]}\mathbf{Dg=0}$ are satisfied (similarly as in %
\cite{ma2}, on N--anholonomic vector bundles, and Chapter 1 in \cite{vsgg},
for generalized Finsler--affine spaces). $\square $
\end{proof}

\begin{theorem}
\label{mconnections}For fixed d--metric, \ $\mathbf{g}_{\alpha \beta
}=[g_{ij},h_{ab}],$ and N--connecti\-on, $\mathbf{N}=\{N_{i}^{a}\},$
structures the set of metric d--connections $^{[M]}\mathbf{\Gamma }_{\
\alpha \beta }^{\gamma }=\widehat{\mathbf{\Gamma }}_{\ \alpha \beta
}^{\gamma }+\ ^{[M]}\mathbf{Z}_{\ \alpha \beta }^{\gamma }$\ is defined by
the deformation d--tensors
\begin{eqnarray*}
^{\lbrack M]}\mathbf{Z}_{\alpha \beta }^{\gamma } &=&\{\ ^{[M]}Z_{\
jk}^{i}=\ ^{[-]}O_{km}^{li}Y_{lj}^{m},\ ^{[M]}Z_{\ bk}^{a}=\
^{[-]}O_{bd}^{ea}Y_{ej}^{m},\  \\
&&\ ^{[M]}Z_{\ ja}^{i}=\ ^{[+]}O_{jk}^{mi}Y_{mc}^{k},\ ^{[M]}Z_{\ bc}^{a}=\
^{[+]}O_{bd}^{ea}Y_{ec}^{d}\}
\end{eqnarray*}%
where the so--called Obata operators are defined
\begin{equation*}
\ ^{[\pm ]}O_{km}^{li}=\frac{1}{2}\left( \delta _{k}^{l}\delta _{m}^{i}\pm
g_{km}g^{li}\right) \mbox{ and }\ ^{[\pm ]}O_{bd}^{ea}=\frac{1}{2}\left(
\delta _{b}^{e}\delta _{d}^{a}\pm h_{bd}h^{ea}\right)
\end{equation*}%
and \ $Y_{lj}^{m},$\ $Y_{ej}^{m},Y_{mc}^{k},$\ $Y_{ec}^{d}$ are arbitrary
d--tensor fields.\
\end{theorem}

\begin{proof}
It also consists from a straightforward verification. Here we note, that \ $%
^{[M]}\mathbf{\Gamma }_{\ \alpha \beta }^{\gamma }$ are generated with
prescribed nontrivial torsion coefficients. If $\ ^{[M]}\mathbf{Z}_{\ \alpha
\beta }^{\gamma }=0,$ the canonical d--connection $\widehat{\mathbf{\Gamma }}%
_{\ \alpha \beta }^{\gamma }$ contains a nonholonomically induced torsion. $%
\square $
\end{proof}

\vskip3pt

We can generalize the concept of N--anholonomic Riemann--Cartan manifold $\
^{RC}\mathbf{V}$ (see Definition \ref{drcm}):

\begin{definition}
\label{dmam}A N--anholonomic metric--affine manifold $^{ma}\mathbf{V}$ is
defined by three fundamental geometric objects: 1) a d--metric $\mathbf{g}%
_{\alpha \beta }=[g_{ij},h_{ab}],$ 2) a N--connection $\mathbf{N}%
=\{N_{i}^{a}\}$ and 3) a general d--connection $\mathbf{D,}$ with nontrivial
nonmetricity d--tensor field $\mathbf{Q=Dg.}$
\end{definition}

The geometry and classification of metric--affine manifolds and related
generalized Finsler--affine spaces is considered in Part I of monograph \cite%
{vsgg}. From Theorems \ref{kmp}, \ref{mconnections} and \ref{teq}, one
follows

\begin{conclusion}
\label{cgmam}The geometry of any manifold $^{ma}\mathbf{V}$ can be
equivalently modelled by deformation tensors on Riemann manifolds provided
with preferred frame structure. The constructions are elaborated in
N--adapted form if we work with the canonical d--connection, or not adapted
to the N--connection structure if we apply the Levi Civita connection.
\end{conclusion}

Finally, in this section, we note that if the torsion and nonmetricity
fields of $^{ma}\mathbf{V}$ are defined by the d--metric and N--connection
coefficients (for instance, in Finsler geometry with Chern or Berwald
connection, see below section \ref{ssxrl}) we can equivalently
(nonholonomically) transform $^{ma}\mathbf{V}$ into a Riemann manifold with
metric structure of type (\ref{metr}) and (\ref{ansatz}).

\section{Einstein Gravity and Lagrange--Finsler Ge\-ometry}

We study N--anholonomic structures in Riemmann--Finsler and Lagrange
geometry modelled on nonholonomic Riemann--Cartan manifolds.

\subsection{Generalized Lagrange spaces}

\label{ssgls}If a N--anholonomic manifold is stated to be a tangent bundle, $%
\mathbf{V=TM},$ the dimension of the base and fiber space coincide, $n=m,$
and we obtain a special case of N--connection geometry \cite{ma1,ma2}. For
such geometric models, a N--connection is defined by Withney sum
\begin{equation}
T\mathbf{TM}=h\mathbf{TM}\oplus v\mathbf{TM,}  \label{nctb}
\end{equation}%
with local coefficients $\mathbf{N}=\{N_{i}^{a}(x^{i},y^{a})\},$ where it is
convenient to distinguish h--indices $i,j,k...$ from v--indices $a,b,c,...$%
\footnote{%
It should be emphasized, that on $\mathbf{TM}$ we can contract h- and
v--indices, which is not possible on a vector bundle $\mathbf{E}$ with $%
n\neq m.$} On $\mathbf{TM,}$ there is an almost complex structure $\mathbf{%
F=\{F}_{\alpha }^{~\beta }\mathbf{\}}$ associated to $\mathbf{N}$ defined by%
\begin{equation}
\mathbf{F(e}_{i}\mathbf{)=-}e_{i}\mbox{\ and }\mathbf{F(}e_{i}\mathbf{)=e}%
_{i},  \label{acs}
\end{equation}%
where $\mathbf{e}_{i}=\partial /\partial x^{i}-N_{i}^{k}\partial /\partial
y^{k}$ and $e_{i}=\partial /\partial y^{i}$ and $\mathbf{F}_{\alpha
}^{~\beta }\mathbf{F}_{\beta }^{~\gamma }=-\delta _{\alpha }^{\beta }.$
Similar constructions can be performed on N--anholonomic manifolds $\mathbf{V%
}^{n+n}$ where fibred structures of dimension $n+n$ are modelled.

A general d--metric structure (\ref{m1}) on $\mathbf{V}^{n+n},$ together
with a prescribed N--connection $\mathbf{N,}$ defines a N--anholonomic
Riemann--Cartan manifold of even dimension.

\begin{definition}
A generalized Lagrange space is modelled on $\mathbf{V}^{n+n}$ (on $\mathbf{%
TM,}$ see \cite{ma1,ma2}) by a d--metric with $g_{ij}=\delta _{i}^{a}\delta
_{j}^{b}h_{ab},$ i.e.%
\begin{equation}
~^{c}\mathbf{g}=\ h_{ij}(x,y)\ \left( e^{i}\otimes e^{j}+\ \ \mathbf{e}%
^{i}\otimes \mathbf{e}^{j}\right) .  \label{sdmgls}
\end{equation}
\end{definition}

One calls $\varepsilon =h_{ab}(x,y)\ y^{a}y^{b}$ to be the absolute energy
associated to a $h_{ab}$ of constant signature.

\begin{theorem}
\label{tcsl1}For nondegenerated Hessians
\begin{equation}
\widetilde{h}_{ab}=\frac{1}{2}\frac{\partial ^{2}\varepsilon }{\partial
y^{a}\partial y^{b}},  \label{hess}
\end{equation}%
when $\det |\widetilde{h}|\neq 0,$ there is a canonical N--connection
completely defined by $\ h_{ij},$%
\begin{equation}
~^{c}N_{i}^{a}(x,y)=\frac{\partial G^{a}}{\partial y^{i}}  \label{cnc}
\end{equation}%
where%
\begin{equation*}
G^{a}=\frac{1}{2}\widetilde{h}^{ab}\left( y^{k}\frac{\partial
^{2}\varepsilon }{\partial y^{b}\partial x^{k}}-\delta _{b}^{k}\frac{%
\partial \varepsilon }{\partial x^{k}}\right) .
\end{equation*}
\end{theorem}

\begin{proof}
One has to consider local coordinate transformation laws for some
coefficients $N_{i}^{a}$ preserving splitting (\ref{nctb}). We can verify
that $~^{c}N_{i}^{a}$ satisfy such conditions. The sketch of proof is given
in \cite{ma1,ma2} for $\mathbf{TM.}$ We can consider any nondegenerated
quadratic form $h_{a^{\prime }b^{\prime }}(x,y)=e_{a^{\prime
}}^{~a}e_{b^{\prime }}^{~b}h_{ab}(x,y)$ on $\mathbf{V}^{n+n}$ if we redefine
the v--coordinates in the form $y^{a^{\prime }}=y^{a^{\prime }}(x^{i},y^{a})$
and $x^{i^{\prime }}=x^{i}.\square $
\end{proof}

Finally, in this section, we state:

\begin{theorem}
\label{tcslg}For any generalized Lagrange space, there are canonical
N--connection $\ ^{c}\mathbf{N,}$ almost complex $^{c}\mathbf{F,}$ d--metric
$~^{c}\mathbf{g}$ and d--connection $~^{c}\widehat{\mathbf{D}}$ structures
defined by an effective regular Lagrangian $~^{\varepsilon }L(x,y)=\sqrt{%
|\varepsilon |}$ and its Hessian $\widetilde{h}_{ab}(x,y)$ (\ref{hess}).
\end{theorem}

\begin{proof}
It follows from formulas (\ref{hess}), (\ref{cnc}), (\ref{acs}) and (\ref%
{cncl}) and adapted d--connection (\ref{cdctb}) and d--metric structures (%
\ref{slm}) all induced by a $~^{\varepsilon }L=\sqrt{|\varepsilon |}.$ $%
\square $
\end{proof}

\vskip3pt

\subsection{Lagrange--Finsler spaces}

The class of Lagrange--Finsler geometries is usually defined on tangent
bundles but it is possible to model such structures on general
N--anholonomic manifolds, for instance, in (pseudo) Riemannian and
Riemann--Cartan geometry, if nonholonomic frames are introduced into
consideration \cite{vsgg,vncg}. Let us consider two such important examples
when the N--anholonomic structures are modelled on $\mathbf{TM.}$ One
denotes by $\widetilde{TM}=TM\backslash \{0\}$ where $\{0\}$ means the set
of null sections of surjective map $\pi :TM\rightarrow M.$

\begin{example}
A Lagrange space is a pair $L^{n}=\left[ M,L(x,y)\right] $ with a
differentiable fundamental Lagrange function $L(x,y)$ defined by a map $%
L:(x,y)\in TM\rightarrow L(x,y)\in \mathbb{R}$ of class $\mathcal{C}^{\infty
}$ on $\widetilde{TM}$ and continuous on the null section $0:M\rightarrow TM$
of $\pi .$ The Hessian (\ref{hess}) is defined
\begin{equation}
\ ^{L}g_{ij}(x,y)=\frac{1}{2}\frac{\partial ^{2}L(x,y)}{\partial
y^{i}\partial y^{j}}  \label{lqf}
\end{equation}%
when $rank\left| g_{ij}\right| =n$ on $\widetilde{TM}$ and the left up ''L''
is an abstract label pointing that certain values are defined by the
Lagrangian $L.$
\end{example}

The notion of Lagrange space was introduced by J. Kern \cite{kern} and
elaborated in details in Ref. \cite{ma1,ma2} as a natural extension of
Finsler geometry. In a more particular case, we have

\begin{example}
\label{exfs}A Finsler space defined by a fundamental Finsler function $%
F(x,y),$ being homogeneous of type $F(x,\lambda y)=\lambda F(x,y),$ for
nonzero $\lambda \in \mathbb{R},$ may be considered as a particular case of
Lagrange geometry when $L=F^{2}.$
\end{example}

Our approach to the geometry of N--anholonomic spaces (in particular, to
that of Lagrange, or Finsler, spaces) is based on canonical d--connections.
It is more related to the existing standard models of gravity and field
theory allowing to define Finsler generalizations of spinor fields,
noncommutative and supersymmetric models, see discussions in \cite{vsgg,vncg}%
. Nevertheless, a number of schools and authors on Finsler geometry prefer
linear connections which are not metric compatible (for instance, the
Berwald and Chern connections, see below Definition \ref{dcbhg}) which
define new classes of geometric models and alternative physical theories
with nonmetricity field, see details in \cite{bcs,bejancu,ma1,ma2,cartan}.
From geometrical point of view, all such approaches are equivalent. It can
be considered as a particular realization, for nonholonomic manifolds, of
the Poincare's idea on duality of geometry and physical models stating that
physical theories can be defined equivalently on different geometric spaces,
see \cite{poincscd}.

From the Theorem \ref{tcslg}, one follows:

\begin{conclusion}
\label{ceqlfg}Any mechanical system with regular Lagrangian $L(x,y)$ (or any
Finsler geometry with fundamental function $F(x,y))$ can be modelled as a
nonhlonomic Riemann geometry with canonical structures $\ ^{L}\mathbf{N,}$ $%
\ ^{L}\mathbf{g}$ and $\ ^{L}\widehat{\mathbf{D}}$ (or$\ ^{F}\mathbf{N,}$ $\
^{F}\mathbf{g}$ and $\ ^{F}\widehat{\mathbf{D}},$ for $L=F^{2})$ defined on
a N--anholonomic manifold $\mathbf{V}^{n+n}.$ In equivalent form, such
Lagrange--Finsler geometries can be described by the same metric and
N--anholonomic distributions but with the corresponding not adapted Levi
Civita connections.
\end{conclusion}

Let us denote by $\mathbf{Ric}(\mathbf{D})=C(1,4)\mathbf{R}(\mathbf{D}),$
where $C(1,4)$ means the contraction on the first and forth indices of the
curvature $\mathbf{R}(\mathbf{D}),$ and $\mathbf{Sc}(\mathbf{D})=C(1,2)%
\mathbf{Ric}(\mathbf{D})=\ ^{s}\mathbf{R},$ where $C(1,2)$ is defined by
contracting $\mathbf{Ric}(\mathbf{D})$ with the inverse d--metric,
respectively, the Ricci tensor and the curvature scalar defined by any
metric d--connection $\mathbf{D}$ and d--metric $\mathbf{g}$ on $\ ^{RC}%
\mathbf{V,}$ see also the component formulas (\ref{dricci}), (\ref{sdccurv})
\ and (\ref{enstdt})\ in Appendix. The Einstein equations are
\begin{equation}
\mathbf{En}(\mathbf{D})\doteqdot \mathbf{Ric}(\mathbf{D})-\frac{1}{2}\mathbf{%
g~Sc}(\mathbf{D})=\mathbf{\Upsilon ,}  \label{einsta}
\end{equation}%
where the source $\mathbf{\Upsilon }$ reflects any contributions of matter
fields and corrections from, for instance, string/brane theories of gravity.
In a physical model, the equations (\ref{einsta}) have to be completed with
equations for the matter fields and torsion (for instance, in the\
Einstein--Cartan theory one considers algebraic equations for the torsion
and its source). It should be noted here that because of nonholonomic
structure of $^{RC}\mathbf{V,}$ the tensor $\mathbf{Ric}(\mathbf{D})$ is not
symmetric and $\mathbf{D}\left[ \mathbf{En}(\mathbf{D})\right] \neq 0.$ This
imposes a more sophisticate form of conservation laws on such spaces with
generic ''local anisotropy'', see discussion in \cite{vsgg} (a similar
situation arises in Lagrange mechanics when nonholonomic constraints modify
the definition of conservation laws). For $\mathbf{D=}\widehat{\mathbf{D}}%
\mathbf{,}$ all constructions can be equivalently redefined for the Levi
Civita connection $\nabla ,$ when $\nabla \left[ En(\nabla )\right] =0.$ A
very important class of models can be elaborated when $\mathbf{\Upsilon =}%
diag\left[ \lambda ^{h}(\mathbf{u})~h\mathbf{g},\lambda ^{v}(\mathbf{u})%
\mathbf{\ }v\mathbf{g}\right] ,$ which defines the so--called N--anholonomic
Einstein spaces with ''nonhomogeneous'' cosmological constant (various
classes of exact solutions in gravity and nonholonomic Ricci flow theory
were constructed and analyzed in \cite{vncg,vsgg,vrf,vv1,vv2}).

\section{Anholonomic Ricci Flows}

The Ricci flow theory was elaborated by R. Hamilton \cite{ham1,ham2}\ and
applied as a method approaching the Poincar\'{e} Conjecture \ and Thurston
Geometrization Conjecture \cite{thur01,thur02}, see Grisha Perelman's works %
\cite{per1,per2,per3} and reviews of results in Refs. \cite{rbook,cao}.

\subsection{Holonomic Ricci flows}

\label{sshrf}

For a one parameter $\tau $ family of Riemannian metrics $\underline{g}(\tau
)=\{\underline{g}_{\alpha \beta }(\tau ,u^{\gamma })\}$ on a N--anholonomic
manifold $\mathbf{V,}$ \ one introduces the Ricci flow equation%
\begin{equation}
\frac{\partial \underline{g}_{\alpha \beta }}{\partial \tau }=-2~_{\shortmid
}\underline{R}_{\alpha \beta },  \label{rfleq}
\end{equation}%
where $_{\shortmid }\underline{R}_{\alpha \beta }$ is the Ricci tensor for
the Levi Civita connection $\bigtriangledown =\{\ _{\shortmid }\Gamma
_{\beta \gamma }^{\alpha }\}$ with the coefficients defined with respect to
a coordinate basis $\partial _{\underline{\alpha }}=\partial /\partial u^{%
\underline{\alpha }}.$ The equation (\ref{rfleq}) \ is a tensor nonlinear
generalization of the scalar heat equation $\partial \phi /\partial \tau
=\bigtriangleup \phi ,$ where $\bigtriangleup $ is the Laplace operator
defined by $\underline{g}.$ Usually, one considers normalized Ricci flows
defined by
\begin{eqnarray}
\frac{\partial }{\partial \tau }g_{\underline{\alpha }\underline{\beta }}
&=&-2\ _{\shortmid }R_{\underline{\alpha }\underline{\beta }}+\frac{2r}{5}g_{%
\underline{\alpha }\underline{\beta }},  \label{feq} \\
g_{\underline{\alpha }\underline{\beta }\mid _{\tau =0}} &=&g_{\underline{%
\alpha }\underline{\beta }}^{[0]}(u),  \label{bc}
\end{eqnarray}%
where the normalizing factor $r=\int RdV/dV$ is introduced in order to
preserve the volume $V,$ the boundary conditions are stated for $\tau =0$
and the solutions are searched for $\tau _{0}>\tau \geq 0.$ For simplicity,
we shall work with equations (\ref{rfleq}) if the constructions will not
result in ambiguities.

It is important to study the evolution of tensors in orthonormal frames and
coframes on nonholonomic manifolds. Let $(\mathbf{V,}g_{\underline{\alpha }%
\underline{\beta }}(\tau )),0\leq \tau <$ $\tau _{0},$ be a Ricci flow with $%
_{\shortmid }\underline{R}_{\alpha \beta }=~$ $_{\shortmid }R_{\underline{%
\alpha }\underline{\beta }}$ and consider the evolution of basis vector
fields
\begin{equation*}
e_{\alpha }(\tau )=e_{\alpha }^{\ \underline{\alpha }}(\tau )~\partial _{%
\underline{\alpha }}\mbox{ and }e_{\ }^{\beta }(\tau )=e_{\ \underline{\beta
}}^{\beta }(\tau )~du^{\underline{\beta }}
\end{equation*}%
which are $\underline{g}(0)$--orthonormal on an open subset $\mathbf{%
U\subset V.}$ We evolve this local frame flows according the formula%
\begin{equation}
\frac{\partial }{\partial \tau }e_{\alpha }^{\ \underline{\alpha }}=g^{%
\underline{\alpha }\underline{\beta }}~_{\shortmid }R_{\underline{\beta }%
\underline{\gamma }}~e_{\alpha }^{\ \underline{\gamma }}.  \label{eof}
\end{equation}%
There are unique solutions for such linear ordinary differential equations
for all time $\tau \in \lbrack 0,\tau _{0}).$

Using the equations (\ref{feq}), (\ref{bc})\ and (\ref{eof}), one can be
defined the evolution equations under Ricci flow, for instance, for the
Riemann tensor, Ricci tensor, Ricci scalar and volume form stated in
coordinate frames (see, for example, the Theorem 3.13 in Ref. \cite{rbook}).
In this section, we shall consider such nonholnomic constrains on the
evolution equation when the geometrical object will evolve in N--adapted
form; we shall also model sets of N--anholnomic geometries, in particular,
flows of geometric objects on nonholonomic Riemann manifolds and Finsler and
Lagrange spaces.

\subsection{Ricci flows and N--anholonomic distributions}

On manifold $\mathbf{V,}$ the equations (\ref{feq}) and (\ref{bc})\ describe
flows not adapted to the N--connections $N_{i}^{a}(\tau ,u).$ For a
prescribed family of such N--connections, we can construct from $\underline{g%
}_{\alpha \beta }(\tau ,u^{\gamma })$ the corresponding set of d--metrics $%
\mathbf{g}_{\alpha \beta }(\tau ,u)$ $\ =[g_{ij}(\tau ,u),h_{ab}(\tau ,u)]$
and the set of N--adapted frames on $(\mathbf{V,g}_{\alpha \beta }(\tau
)),0\leq \tau <$ $\tau _{0}.$ The evolution of such N--adapted frames is
defined not by the equations (\ref{eof}) but satisfies the

\begin{proposition}
For a prescribed $n+m$ splitting, the solutions of the system (\ref{feq})
and (\ref{bc}) define a natural flow of preferred N--adapted frame
structures.
\end{proposition}

\begin{proof}
Following formulas (\ref{metr}), (\ref{ansatz}) and (\ref{aux1a}), the
boundary conditions (\ref{bc}) state the values $N_{i}^{a}(\tau =0,u)$ and $%
\mathbf{g}_{\alpha \beta }(\tau =0,u)=[g_{ij}(\tau =0,u),h_{ab}(\tau =0,u)].$
Having a well defined solution $\underline{g}_{\alpha \beta }(\tau ,u),$ we
can construct the coefficients of N--connection $N_{i}^{a}(\tau ,u)$ and
d--metric $\mathbf{g}(\tau ,u)$ $\ =[g(\tau ,u),h(\tau ,u)]$ for any $\tau
\in \lbrack 0,\tau _{0}):$ \ the associated set of frame (vielbein)
structures $\mathbf{e}_{\nu }(\tau )=(\mathbf{e}_{i}(\tau ),e_{a}),$ where
\begin{equation}
\mathbf{e}_{i}(\tau )=\frac{\partial }{\partial x^{i}}-N_{i}^{a}(\tau ,u)%
\frac{\partial }{\partial y^{a}}\mbox{ and
}e_{a}=\frac{\partial }{\partial y^{a}},  \label{dder1s}
\end{equation}%
and the set of dual frame (coframe) structures $\mathbf{e}^{\mu }(\tau
)=(e^{i},\mathbf{e}^{a}(\tau )),$ where
\begin{equation}
e^{i}=dx^{i}\mbox{ and }\mathbf{e}^{a}(\tau )=dy^{a}+N_{i}^{a}(\tau
,u)dx^{i}.  \label{ddif1s}
\end{equation}%
$\square $
\end{proof}

\vskip3pt

We conclude that prescribing the existence of a nonintegrable $(n+m)$%
--decomposition on a manifold for any $\tau \in \lbrack 0,\tau _{0}),$ from
any solution of the Ricci flow equations (\ref{eof}), we can extract a set
of preferred frame structures with associated N--connections, with respect
to which we can perform the geometric constructions in N--adapted form.

We shall need a forumula relating the connection Laplacian on contravariant
one--tensors with Ricci curvature and the corresponding deformations under
N--anholonomic maps. Let $\mathbf{A}$ be a d--tensor of rank $k.$ Then we
define $\nabla ^{2}\mathbf{A,}$ for $\nabla $ being the Levi Civita
connection, to be a contravariant tensor of rank $k+2$ given by
\begin{equation}
\nabla ^{2}\mathbf{A(\cdot ,X,Y)=(\nabla }_{\mathbf{X}}\mathbf{\nabla }_{%
\mathbf{Y}}\mathbf{A)(\cdot )-}\nabla _{\mathbf{\nabla }_{\mathbf{X}}\mathbf{%
Y}}\mathbf{A(\cdot ).}  \label{auxil1a}
\end{equation}%
This defines the (Levi Civita) connection Laplacian
\begin{equation}
\Delta \mathbf{A\doteqdot g}^{\alpha \beta }\left( \nabla ^{2}\mathbf{A}%
\right) \left( \mathbf{e}_{\alpha },\mathbf{e}_{\beta }\right) ,
\label{auxil1d}
\end{equation}%
for tensors, and
\begin{equation*}
\Delta f\mathbf{\doteqdot }tr~\nabla ^{2}f=\mathbf{g}^{\alpha \beta }\left(
\nabla ^{2}f\right) _{\alpha \beta },
\end{equation*}%
for a scalar function on $\mathbf{V.}$ In a similar manner, by substituting $%
\nabla $ with $\widehat{\mathbf{D}},$ we can introduce the canonical
d--connection Laplacian, for instance,
\begin{equation}
\widehat{\Delta }\mathbf{A\doteqdot g}^{\alpha \beta }\left( \widehat{%
\mathbf{D}}^{2}\mathbf{A}\right) \left( \mathbf{e}_{\alpha },\mathbf{e}%
_{\beta }\right) .  \label{auxil1b}
\end{equation}

\begin{proposition}
\label{propdefl}The Laplacians $\widehat{\Delta }$ and $\Delta $ are related
by formula%
\begin{equation}
\Delta \mathbf{A}=\widehat{\Delta }\mathbf{A}+~_{\shortmid }\Delta \mathbf{A}
\label{auxil1c}
\end{equation}%
where the deformation d--tensor of the Laplacian, $~_{\shortmid }\Delta ,$
is defined canonically by the N--connection and d--metric coefficients.
\end{proposition}

\begin{proof}
We sketch the method of computation $~_{\shortmid }\Delta .$ Using the
formula (\ref{cdeft}), we have%
\begin{equation}
\mathbf{\nabla }_{\mathbf{X}}=\widehat{\mathbf{D}}_{\mathbf{X}}+\
_{\shortmid }Z_{\mathbf{X}}  \label{auxil1}
\end{equation}%
where$\ \ _{\shortmid }Z_{\mathbf{X}}=\mathbf{X}^{\alpha }~\ _{\shortmid
}Z_{\ \alpha \beta }^{\gamma }$ is defined for any $\mathbf{X}^{\alpha }$
with $\ _{\shortmid }\mathbf{Z}_{\ \alpha \beta }^{\gamma }$ computed
following formulas (\ref{cdeft}); all such coefficients depend on
N--connection and d--metric coefficients and their derivatives, i.e. on
generic off--diagonal metric coefficients (\ref{ansatz}) and their
derivatives. Introducing (\ref{auxil1}) into (\ref{auxil1a}) and (\ref%
{auxil1d}), and separating the terms depending only on $\widehat{\mathbf{D}}%
_{\mathbf{X}}$ we get $\widehat{\Delta }\mathbf{A}$ (\ref{auxil1b}). The
rest of terms with linear or quadratic dependence on $\ _{\shortmid }Z_{\
\alpha \beta }^{\gamma }$ and their derivatives define
\begin{equation*}
~_{\shortmid }\Delta \mathbf{A\doteqdot g}^{\alpha \beta }\left( ~_{\Delta }%
\mathbf{ZA}\right) ,
\end{equation*}%
where
\begin{eqnarray*}
~_{\Delta }\mathbf{ZA} &\mathbf{=}&\widehat{\mathbf{D}}_{\mathbf{X}}\left( \
_{\shortmid }Z_{\mathbf{Y}}\mathbf{A}\right) +\ _{\shortmid }Z_{\mathbf{X}%
}\left( \widehat{\mathbf{D}}_{\mathbf{Y}}\mathbf{A}\right) +\ _{\shortmid
}Z_{\mathbf{X}}\left( \ _{\shortmid }Z_{\mathbf{Y}}\mathbf{A}\right) \\
&&-\widehat{\mathbf{D}}_{\ _{\shortmid }Z_{\mathbf{Y}}}\mathbf{A-}\
_{\shortmid }Z_{\widehat{\mathbf{D}}_{\mathbf{X}}\mathbf{Y}}\mathbf{A-}\
_{\shortmid }Z_{\ _{\shortmid }Z_{\mathbf{X}}\mathbf{Y}}\mathbf{A.}
\end{eqnarray*}%
$\square $
\end{proof}

\vskip3pt

In a similar form as for Proposition \ref{propdefl}, we prove

\begin{proposition}
\label{propdefla}The curvature, Ricci and scalar tensors of the Levi Civita
connection $\nabla $ and the canonical d--connection $\widehat{\mathbf{D}}$
are defined by formulas%
\begin{eqnarray*}
~_{\shortmid }R\mathbf{(X,Y)} &=&\widehat{\mathbf{R}}\mathbf{(X,Y)+}%
~_{\shortmid }\widehat{Z}\mathbf{(X,Y),} \\
Ric(\mathbf{\nabla }) &=&\mathbf{Ric}(\widehat{\mathbf{D}})+Ric(\mathbf{%
~_{\shortmid }}\widehat{Z}), \\
Sc(\mathbf{\nabla }) &=&\mathbf{Sc}(\widehat{\mathbf{D}})+Sc(\mathbf{%
~_{\shortmid }}\widehat{Z}),
\end{eqnarray*}%
where%
\begin{eqnarray*}
\mathbf{~_{\shortmid }}\widehat{Z}\mathbf{(X,Y)} &=&\mathbf{D}_{\mathbf{X}}\
_{\shortmid }Z_{\mathbf{Y}}-\ _{\shortmid }Z_{\mathbf{Y}}\mathbf{D}_{\mathbf{%
X}}\mathbf{-}\ _{\shortmid }Z_{\mathbf{[X,Y]}} \\
Ric(\mathbf{~_{\shortmid }}\widehat{\mathbf{Z}}) &=&C(1,4)\mathbf{%
~_{\shortmid }}\widehat{Z},~Sc(\mathbf{~_{\shortmid }}\widehat{\mathbf{Z}}%
)=C(1,2)Ric(\mathbf{~_{\shortmid }}\widehat{Z})
\end{eqnarray*}%
for $\widehat{\mathbf{R}}$ computed following formula (\ref{curv1}) and $%
\mathbf{Sc}(\widehat{\mathbf{D}})=\ ^{s}\widehat{\mathbf{R}}$.
\end{proposition}

In the theory of Ricci flows, one consider tensors quadratic in the
curvature tensors, for instance, for any given $\mathbf{g}^{\beta \beta
^{\prime }}$ and $\mathbf{D}$
\begin{eqnarray}
\mathbf{B}_{\alpha \gamma \alpha ^{\prime }\gamma ^{\prime }} &=&\mathbf{g}%
^{\beta \beta ^{\prime }}\mathbf{g}^{\delta \delta ^{\prime }}\mathbf{R}%
_{\alpha \beta \gamma \delta }\mathbf{R}_{\alpha ^{\prime }\beta ^{\prime
}\gamma ^{\prime }\delta ^{\prime }},  \label{aux5} \\
\underline{\mathbf{B}}_{\alpha \gamma \alpha ^{\prime }\gamma ^{\prime }}
&\doteqdot &\mathbf{B}_{\alpha \gamma \alpha ^{\prime }\gamma ^{\prime }}-%
\mathbf{B}_{\alpha \gamma \alpha ^{\prime }\gamma ^{\prime }}-\mathbf{B}%
_{\alpha \gamma ^{\prime }\gamma \alpha ^{\prime }}+\mathbf{B}_{\alpha
\gamma \alpha ^{\prime }\gamma ^{\prime }},  \notag \\
\underline{\mathbf{B}}_{\alpha \gamma ^{\prime }} &\doteqdot &\mathbf{D}%
_{\alpha }\mathbf{D}_{\gamma ^{\prime }}\ ^{s}\mathbf{R-g}^{\beta \beta
^{\prime }}\left( \mathbf{D}_{\beta ^{\prime }}\mathbf{D}_{\alpha }\mathbf{R}%
_{\gamma ^{\prime }\beta }+\mathbf{D}_{\beta ^{\prime }}\mathbf{D}_{\gamma
^{\prime }}\mathbf{R}_{\alpha \beta }\right) .  \notag
\end{eqnarray}%
Using the connections $\nabla ,$ or $\widehat{\mathbf{D}}\mathbf{,}$ we
similarly define and compute the values $\mathbf{~_{\shortmid }}B_{\alpha
\gamma \alpha ^{\prime }\gamma ^{\prime }},\mathbf{~_{\shortmid }}\underline{%
B}_{\alpha \gamma \alpha ^{\prime }\gamma ^{\prime }}$ and $\mathbf{%
~_{\shortmid }}\underline{B}_{\alpha \gamma ^{\prime }},$ or $\widehat{%
\mathbf{B}}_{\alpha \gamma \alpha ^{\prime }\gamma ^{\prime }},\widehat{%
\underline{\mathbf{B}}}_{\alpha \gamma \alpha ^{\prime }\gamma ^{\prime }}$
and $\widehat{\underline{\mathbf{B}}}_{\alpha \gamma ^{\prime }}.$

\subsection{Evolution of distinguished geometric objects}

There are d--objects (d--tensors, d--connections) with N--adapted evolution
completely defined by solutions of the Ricci flow equations (\ref{eof}).

\begin{definition}
A geometric structure/object is extracted from a (Riemanni\-an) Ricci flow
(for the Levi Civita connection) if the corresponding structure/object can
be redefined equivalently, prescribing a $(n+m)$--splitting, as a N--adapted
structure/ d--object subjected to corresponding N--anholonomic flows.
\end{definition}

Following the Propositions \ref{propdefl} and \ref{propdefla} (we emphasize
the calculus used as proofs) and formulas (\ref{aux5}), we prove

\begin{theorem}
\label{th313}The evolution equations for the Riemann and Ricci tensors and
scalar curvature defined by the canonical d--connection are extracted
respectively:
\begin{eqnarray*}
\frac{\partial }{\partial \tau }\widehat{\mathbf{R}}_{\alpha \beta \gamma
\delta } &=&\widehat{\Delta }\widehat{\mathbf{R}}_{\alpha \beta \gamma
\delta }+2\widehat{\underline{\mathbf{B}}}_{\alpha \beta \gamma \delta }+%
\widehat{\mathbf{Q}}_{\alpha \beta \gamma \delta }, \\
\frac{\partial }{\partial \tau }\widehat{\mathbf{R}}_{\alpha \beta } &=&%
\widehat{\Delta }\widehat{\mathbf{R}}_{\alpha \beta }+\widehat{\mathbf{Q}}%
_{\alpha \beta }, \\
\frac{\partial }{\partial \tau }\ ^{s}\widehat{\mathbf{R}} &=&\widehat{%
\Delta }\ ^{s}\widehat{\mathbf{R}}+2\widehat{\mathbf{R}}_{\alpha \beta }%
\widehat{\mathbf{R}}^{\alpha \beta }+\widehat{\mathbf{Q}}
\end{eqnarray*}%
where, for%
\begin{eqnarray*}
\mathbf{_{\shortmid }}R_{\alpha \beta \gamma \delta } &=&\widehat{\mathbf{R}}%
_{\alpha \beta \gamma \delta }+\mathbf{~_{\shortmid }}Z_{\alpha \beta \gamma
\delta },~\underline{B}_{\alpha \beta \gamma \delta }=\widehat{\underline{%
\mathbf{B}}}_{\alpha \beta \gamma \delta }+\mathbf{~_{\shortmid }}\underline{%
\widehat{Z}}_{\alpha \beta \gamma \delta },~Z=\mathbf{g}^{\alpha \beta }%
\mathbf{~_{\shortmid }}Z_{\alpha \beta }, \\
\mathbf{_{\shortmid }}R_{\alpha \beta } &=&\widehat{\mathbf{R}}_{\alpha
\beta }+\mathbf{~_{\shortmid }}Z_{\alpha \beta },\mathbf{~_{\shortmid }}%
\underline{B}_{\alpha \gamma ^{\prime }}=\widehat{\underline{\mathbf{B}}}%
_{\alpha \gamma ^{\prime }}+\mathbf{~_{\shortmid }}\underline{\widehat{Z}}%
_{\alpha \gamma ^{\prime }},\ \mathbf{~_{\shortmid }}^{s}R=\ ^{s}\widehat{%
\mathbf{R}}+Z,
\end{eqnarray*}%
the $Q$--terms (defined by the coefficients of canonical d--connection, $%
N_{i}^{a}$ and $\mathbf{g}_{\alpha \beta }=\left[ g_{ij},h_{ab}\right] $ and
their derivatives) are
\begin{eqnarray*}
\widehat{\mathbf{Q}}_{\alpha \beta \gamma \delta } &=&-\frac{\partial }{%
\partial \tau }\mathbf{~_{\shortmid }}Z_{\alpha \beta \gamma \delta
}+~_{\shortmid }\Delta \widehat{\mathbf{R}}_{\alpha \beta \gamma \delta }+2%
\mathbf{~_{\shortmid }}\underline{\widehat{Z}}_{\alpha \beta \gamma \delta },
\\
\widehat{\mathbf{Q}}_{\alpha \beta } &=&-\frac{\partial }{\partial \tau }%
\mathbf{~_{\shortmid }}Z_{\alpha \beta }+~_{\shortmid }\Delta \widehat{%
\mathbf{R}}_{\alpha \beta }+\mathbf{~_{\shortmid }}\underline{\widehat{Z}}%
_{\alpha \beta }, \\
\widehat{\mathbf{Q}} &=&-\frac{\partial }{\partial \tau }Z+\widehat{\Delta }%
Z+~_{\shortmid }\Delta \ ^{s}\widehat{\mathbf{R}}+2\widehat{\mathbf{R}}%
_{\alpha \beta }\mathbf{~_{\shortmid }}Z^{\alpha \beta }+2\mathbf{%
~_{\shortmid }}Z_{\alpha \beta }\widehat{\mathbf{R}}^{\alpha \beta }+2%
\mathbf{~_{\shortmid }}Z_{\alpha \beta }\mathbf{~_{\shortmid }}Z^{\alpha
\beta }.
\end{eqnarray*}
\end{theorem}

In Ricci flow theory, it is important to have the formula for the evolution
of the volume form:

\begin{remark}
\label{th313a}The\ deformation of the volume form is stated by equation
\begin{equation*}
\frac{\partial }{\partial \tau }dvol\left( \tau ,u^{\alpha }\right) =-\left(
\ ^{s}\widehat{\mathbf{R}}+Z\right) \ dvol\left( \tau ,u^{\alpha }\right)
\end{equation*}%
which is just that for the Levi Civita connection and
\begin{equation*}
dvol\left( \tau ,u^{\alpha }\right) \doteqdot \sqrt{|\det \underline{g}%
_{\alpha \beta }\left( \tau ,u^{\gamma }\right) |},
\end{equation*}%
where.$\underline{g}_{\alpha \beta }(\tau )$ are metrics of type (\ref{metr}%
).
\end{remark}

The evolution equations from Theorem \ref{th313} and Remark \ref{th313a}
transform into similar ones from Theorem 3.13 in Ref. \cite{rbook}.

For any solution of equations (\ref{feq}) and (\ref{bc}), on $\mathbf{%
U\subset V,}$ we can construct for any $\tau \in \lbrack 0,\tau _{0})$ a
parametrized set of canonical d--connections $\widehat{\mathbf{D}}(\tau )=\{%
\widehat{\mathbf{\Gamma }}_{\ \alpha \beta }^{\gamma }(\tau )\}$ (\ref%
{candcon}) defining the corresponding canonical Riemann d--tensor (\ref%
{dcurv}), nonsymmetric Ricci d--tensor $\widehat{\mathbf{R}}_{\alpha \beta }$
(\ref{dricci}) and scalar (\ref{sdccurv}). The coefficients of d--objects
are defined with respect to evolving N--adapted frames (\ref{dder1s}) and (%
\ref{ddif1s}). One holds

\begin{conclusion}
The evolution of corresponding d--objects on N--anholonom\-ic Riemann
manifolds can be canonically extracted from the evolution under Ricci flows
of geometric objects on Riemann manifolds.
\end{conclusion}

In the sections \ref{ssxrf} and \ref{ssxrl}, we shall consider how Finsler
and Lagrange configurations can be extracted by more special
parametrizations of metric and nonholonomic constraints.

\subsection{Nonholonomic Ricci flows of (non)symmetric metrics}

The Ricci flow equations were introduced by R. Hamilton \cite{ham1} in a
heuristic form similarly to that how A. Einstein proposed his equations by
considering possible physically grounded equalities between the metric and
its first and second derivatives and the second rank Ricci tensor. On
(pseudo) Riemannian spaces the metric and Ricci tensors are both symmetric
and it is possible to consider the parameter derivative of metric and/or
correspondingly symmetrized energy--momentum of matter fields as sources for
the Ricci tensor.

On N--anholonomic manifolds there are two alternative possibilities: The
first one is to postulate the Ricci flow equations in symmetric form, for
the Levi Civita connection, and then to extract various N--anholonomic
configurations by imposing corresponding nonholonomic constraints. The bulk
of our former and present work are related to symmetric metric
configurations.

In the second case, we can start from the very beginning with a nonsymmetric
Ricci tensor for a non--Riemannian space. In this section, we briefly
speculate on such geometric constructions: \ The nonholonomic Ricci flows
even beginning with a symmetric metric tensor may result naturally in
nonsymmetric metric tensors $\widehat{g}_{\alpha \beta }=\underline{g}%
_{\alpha \beta }+\overleftrightarrow{g}_{\alpha \beta },$ where $%
\overleftrightarrow{g}_{\alpha \beta }=-\overleftrightarrow{g}_{\beta \alpha
}.$ Nonsymmetric metrics in gravity were originally considered by A.\
Einstein \cite{einstns} and L. P. Eisenhart \cite{eisenh}, see modern
approaches in Ref. \cite{moffat}. For Finsler and Lagrange spaces, such
nonsymmetric metric generalizations were performed originally in Refs. \cite%
{matan,asm} (Chapter 8 of monograph \cite{ma1} contains a review of results
on Eisenhart--Lagrange spaces).

\begin{theorem}
With respect to N--adapted frames, the canonical nonholonomic Ricci flows
with nonsymmetric metrics defined by equations
\begin{eqnarray}
\frac{\partial }{\partial \tau }g_{ij} &=&-2\widehat{R}_{ij}+2\lambda
g_{ij}-h_{cd}\frac{\partial }{\partial \tau }(N_{i}^{c}N_{j}^{d}),
\label{1eq} \\
\frac{\partial }{\partial \tau }h_{ab} &=&-2\widehat{R}_{ab}+2\lambda
h_{ab},\   \label{2eq} \\
\frac{\partial }{\partial \tau }\overleftrightarrow{g}_{ia} &=&\widehat{R}%
_{ia},~\text{ }\frac{\partial }{\partial \tau }\overleftrightarrow{g}_{ai}=%
\widehat{R}_{ai}  \label{3eq}
\end{eqnarray}%
where $\underline{g}_{\alpha \beta }=[g_{ij},h_{ab}]$ with respect to
N--adapted basis (\ref{ddif}), $\lambda =r/5,$ $y^{3}=v$ and $\tau $ can be,
for instance, the time like coordinate, $\tau =t,$ or any parameter or extra
dimension coordinate.
\end{theorem}

\begin{proof}
It follows from a redefinition of equations (\ref{feq}) with respect to
N--adapted frames (by using the frame transform (\ref{vt1}) and (\ref{vt2}%
)), and considering respectively the canonical Ricci d--tensor (\ref{dricci}%
) constructed from $[g_{ij},h_{ab}].$ Here we note that normalizing factor $%
r $ is considered for the symmetric part of metric. $\square $
\end{proof}

\vskip3pt

One follows:

\begin{conclusion}
Nonholonomic Ricci flows (for the canonical d--connection) resulting in
symmetric d--metrics are parametrized by the constraints
\begin{equation}
\overleftrightarrow{g}_{\alpha \beta }=0\mbox{ and }\widehat{R}_{ia}=%
\widehat{R}_{ai}=0.  \label{3eqs}
\end{equation}
\end{conclusion}

The system of equations (\ref{1eq}), (\ref{2eq}) and (\ref{3eqs}), for
''symmetric'' nonholonomic Ricci flows, was introduced and analyzed in Refs. %
\cite{vrf,vv1}.

\begin{example}
The conditions (\ref{3eqs}) are satisfied by any ansatz of type (\ref{m1})
in 3D, 4D , or 5D, with coefficients of type
\begin{equation}
g_{i}=g_{i}(x^{k}),h_{a}=h_{a}(x^{k},v),N_{i}^{3}=w_{i}(x^{k},v),N_{i}^{4}=n_{i}(x^{k},v),
\label{m2ac}
\end{equation}%
for $i,j,...=1,2,3$ and $a,b,...=4,5$ (the 3D and 4D being parametrized by
eliminating the cases $i=1$ and, respectively, $i=1,2);$ $y^{4}=v$ being the
so--called ''anisotropic'' coordinate. Such metrics are off--diagonal with
the coefficients depending on 2 and 3 coordinates but positively not
depending on the coordinate $y^{5}.$
\end{example}

We constructed and investigated various types of exact solutions of the
nonholonomc Einstein equations and Ricci flow equations, respectively in
Refs. \cite{vncg,vsgg,vcla} and \cite{vrf,vv1,vv2}. They are parametrized by
ansatz of type (\ref{m2ac}) which positively constrains the Ricci flows to
be with symmetric metrics. Such solutions can be used as backgrounds for
investigating flows of Eisenhart (generalized Finsler--Eisenhart geometries)
if the constraints (\ref{3eqs}) are not completely imposed. We shall not
analyze this type of N--anholonomic Ricci flows in this series of works.

\section{Generalized Finsler--Ricci Flows}

The aim of this section is to provide some examples illustrating how
different types of nonholonomic constraints on Ricci flows of Riemannian
metrics model different classes of N--anholonomic spaces (defined by Finsler
metrics and connections, geometric models of Lagrange mechanics and
generalized Lagrange geometries).

\subsection{Finsler--Ricci flows}

\label{ssxrl}Let us consider a $\tau $--parametrized family (set) of
fundamental Finsler functions $F(\tau )=F(\tau ,x^{i},y^{a}),$ see Example %
\ref{exfs}.\footnote{%
we shall write, in brief, only the parameter dependence and even omit
dependencies both on coordinates and parameter if that will not result in
ambiguities} For a family of nondegenerated Hessians

\begin{equation}
\ ^{F}h_{ij}(\tau ,x,y)=\frac{1}{2}\frac{\partial ^{2}F^{2}(\tau ,x,y)}{%
\partial y^{i}\partial y^{j}},  \label{ffh}
\end{equation}%
see formula (\ref{lqf}) for effective $\varepsilon (\tau )=L(\tau
)=F^{2}(\tau ),$ we can model Finsler metrics on $\mathbf{V}^{n+n}$ (or on $%
\mathbf{TM)}$ and the corresponding family of canonical N--connections, see (%
\ref{cnc}),
\begin{equation}
~^{c}N_{i}^{a}(\tau )=\frac{\partial G^{a}(\tau )}{\partial y^{i}},
\label{cncff}
\end{equation}%
where%
\begin{equation*}
G^{a}(\tau )=\frac{1}{2}\ ^{F}h^{ab}(\tau )\left( y^{k}\frac{\partial
^{2}F^{2}(\tau )}{\partial y^{b}\partial x^{k}}-\delta _{b}^{k}\frac{%
\partial F^{2}(\tau )}{\partial x^{k}}\right)
\end{equation*}%
and $\ ^{F}h^{ab}(\tau )$ are inverse to $\ ^{F}h_{ij}(\tau ).$

\begin{proposition}
Any family of fundamental Finsler functions $F(\tau )$ with nondegenerated $%
\ ^{F}h_{ij}(\tau )$ defines a corresponding family of Sasaki type metrics
\begin{equation}
~^{c}\mathbf{g}(\tau )=\ \ ^{F}h_{ij}(\tau ,x,y)\ \left( e^{i}\otimes
e^{j}+\ ~^{c}\mathbf{e}^{i}(\tau )\otimes ~^{c}\mathbf{e}^{j}(\tau )\right) ,
\label{sdmglsff}
\end{equation}%
with $\ ^{F}g_{ij}(\tau )=\ ^{F}h_{ij}(\tau ,x,y),\ $where $\ ~^{c}\mathbf{e}%
^{a}(\tau )=dy^{a}+~^{c}N_{i}^{a}(\tau ,u)dx^{i}$ are defined by the
N--connection (\ref{cncff}).
\end{proposition}

\begin{proof}
It follows from \ the explicit construction (\ref{sdmglsff}).$\square $
\end{proof}

\vskip3pt

For $\mathbf{V}^{n+n}=\mathbf{TM}=(TM,\pi ,M,~^{c}N_{i}^{a})\ $with
injective $\pi :TM\rightarrow M,$ we can model by $F(\tau )$ various classes
of Finsler geometries. In explicit form, we work on $\widetilde{TM}\doteqdot
TM\backslash \{0\}$ and consider the pull--buck bundle $\pi ^{\ast }TM.$ One
generates sets of geometric objects on pull--back cotangent bundle $\pi
^{\ast }T^{\ast }M$ and its tensor products:

\begin{description}
\item on $\pi ^{\ast }T^{\ast }M\otimes \pi ^{\ast }T^{\ast }M\otimes \pi
^{\ast }T^{\ast }M,$ a corresponding family of Cartan tensors%
\begin{eqnarray*}
A(\tau ) &=&A_{ijk}(\tau )dx^{i}\otimes dx^{j}\otimes dx^{k}, \\
A_{ijk}(\tau ) &\doteqdot &\frac{F(\tau )}{2}\frac{\partial g_{ij}(\tau )}{%
\partial y^{k}};
\end{eqnarray*}

\item on $\pi ^{\ast }T^{\ast }M,$ a family of Hilbert forms
\begin{equation*}
\omega (\tau )\doteqdot \frac{\partial F(\tau )}{\partial y^{k}}dx^{i},
\end{equation*}%
and the d--connection 1--form%
\begin{eqnarray}
\omega _{j}^{~i}(\tau ) &=&L_{\ jk}^{i}(\tau )dx^{k}  \label{hform} \\
L_{\ jk}^{i}(\tau ) &=&\frac{1}{2}\ ^{F}g^{ih}(\ ~^{c}\mathbf{e}_{k}\
^{F}g_{jh}+\ ~^{c}\mathbf{e}_{j}\ ^{F}g_{kh}-\ ~^{c}\mathbf{e}_{h}\
^{F}g_{jk}).  \notag
\end{eqnarray}
\end{description}

\begin{theorem}
The set of fundamental Finsler functions $F(\tau )$ defines on $\pi ^{\ast
}TM$ a unique set of linear connections, called the Chern connections,
characterized by the structure equations:%
\begin{equation*}
d(dx^{i})-dx^{i}\wedge \omega _{i}^{~i}(\tau )=-dx^{i}\wedge \omega
_{i}^{~i}(\tau )=0,
\end{equation*}%
i.e. the torsion free condition;%
\begin{equation*}
dg_{ij}(\tau )-\ ^{F}g_{kj}(\tau )\omega _{i}^{~k}(\tau )-\ ^{F}g_{ik}(\tau
)\omega _{j}^{~k}(\tau )=2\frac{A_{ija}(\tau )}{F(\tau )}\ ~^{c}\mathbf{e}%
^{a}(\tau ),
\end{equation*}%
i.e. the almost metric compatibility condition.
\end{theorem}

\begin{proof}
It follows from \ straightforward computations. For any fixed value $\tau
=\tau _{0},$ it is just the Chern's Theorem 2.4.1. from \cite{bcs}. $\square
$
\end{proof}

\vskip3pt

In order to elaborate a complete geometric model on $TM,$ which also allows
us to perform the constructions for N--anholonomic manifolds, we have to
extend the above considered forms with nontrivial coefficients with respect
to $\ ~^{c}\mathbf{e}^{a}(\tau ).$

\begin{definition}
\label{dcbhg}A family of fundamental Finsler metrics $F(\tau )$ defines
models of Finsler geometry (equivalently, space) with d--connections $\Gamma
_{\ \beta \gamma }^{\alpha }(\tau )=(L_{\ jk}^{i}(\tau ),$ $C_{\
jk}^{i}(\tau ))$ on a corresponding N--anholonomic manifold $\mathbf{V:}$

\begin{itemize}
\item of Cartan type if $L_{\ jk}^{i}(\tau )$ is that from (\ref{hform}) and
\begin{equation}
C_{\ jk}^{i}(\tau )=\frac{1}{2}\ ^{F}g^{ih}(\frac{\partial }{\partial x^{k}}%
\ ^{F}g_{jh}+\frac{\partial }{\partial x^{j}}\ ^{F}g_{kh}-\frac{\partial }{%
\partial x^{h}}\ ^{F}g_{jk}),  \label{vcomp}
\end{equation}%
which is similar to formulas (\ref{cdctb}) but for $L=F^{2}(\tau );$

\item of Chern type if $L_{\ jk}^{i}(\tau )$ is given by (\ref{hform}) and $%
C_{\ jk}^{i}(\tau )=0;$

\item of Berwald type if $L_{\ jk}^{i}(\tau )=\partial
~^{c}N_{j}^{i}/\partial y^{k}$ and $C_{\ jk}^{i}(\tau )=0;$

\item of Hashiguchi type if $L_{\ jk}^{i}(\tau )=\partial
~^{c}N_{j}^{i}/\partial y^{k}$ and $C_{\ jk}^{i}(\tau )$ is given by (\ref%
{vcomp}).
\end{itemize}
\end{definition}

Various classes of remarkable Finsler connections have been investigated in
Refs. \cite{ma1,ma2,bejancu,bcs}, see \cite{vncg,vsgg} on modelling Finsler
like structures in Einstein and string gravity and in noncommutative gravity.

It should be emphasized that the models of Finsler geometry with Chern,
Berwald or Hashiguchi type d--connections are with nontrivial nonmetricity
field. So, in general, a family of Finsler fundamental metric functions $%
F(\tau )$ may generate various types of N--anholonomic metric--affine
geometric configurations, see Definition \ref{dmam}, but all components of
such induced nonmetricity and/or torsion fields are defined by the
coefficients of corresponding families of generic off--diagonal metrics of
type (\ref{metr}), when the ansatz (\ref{ansatz}) is parametrized for $%
g_{ij}=h_{ij}=\ ^{F}h_{ij}(\tau )$ and $N_{i}^{a}=$ $~^{c}N_{i}^{a}(\tau ).$
Applying the results of Theorem \ref{mconnections}, we can transform the
families of ''nonmetric'' Finsler geometries into corresponding metric ones
and model the Finsler configurations on N--anholonomic Riemannian spaces,
see Conclusion \ref{cgmam}. In the ''simplest'' geometric and physical
manner (convenient both for applying the former Hamilton--Perelman results
on Ricci flows for Riemannian metrics, as well for further generalizations
to noncommutative Finsler geometry, supersymmetric models and so on...), we
restrict our analysis to Finsler--Ricci flows with canonical d--connection
of Cartan type when $~^{F}\widehat{\mathbf{\Gamma }}_{\ \beta \gamma
}^{\alpha }(\tau )=(L_{\ jk}^{i}(\tau ),$ $C_{\ jk}^{i}(\tau ))$ is with $%
L_{\ jk}^{i}(\tau )$ from (\ref{hform}) and $C_{\ jk}^{i}(\tau )$ from (\ref%
{vcomp}). This provides a proof for

\begin{lemma}
\label{ldconff}A family of Finsler geometries defined by $F(\tau )$ can be
characterized equivalently by the corresponding canonical d--connections (in
N--adapted form) and Levi Civita connections (in not N--adapted form)
related by formulas%
\begin{equation}
\ _{\shortmid }^{F}\Gamma _{\ \alpha \beta }^{\gamma }=~^{F}\widehat{\mathbf{%
\Gamma }}_{\ \alpha \beta }^{\gamma }+\ _{\shortmid }Z_{\ \alpha \beta
}^{\gamma }  \label{fdlc}
\end{equation}%
where $\ _{\shortmid }Z_{\ \alpha \beta }^{\gamma }$ is computed following
formulas (\ref{cdeftc}) for $g_{ij}=h_{ij}=\ ^{F}h_{ij}(\tau )$ and $%
N_{i}^{a}=$ $~^{c}N_{i}^{a}(\tau ).$
\end{lemma}

Following the Lemma \ref{ldconff} and section \ref{sshrf}, we obtain the
proof of

\begin{theorem}
\label{tfrf}The Finsler--Ricci flows for fundamental metric functions\newline
$F(\tau )$ can be extracted from usual Ricci flows of Riemannian metrics
paramet\-riz\-ed in the form%
\begin{equation}
\ ^{F}g_{\underline{\alpha }\underline{\beta }}(\tau )=\left[
\begin{array}{cc}
\ ^{F}g_{ij}+~^{c}N_{i}^{a}~^{c}N_{j}^{b}\ ^{F}g_{ab} & ~^{c}N_{j}^{e}\
^{F}g_{ae} \\
~^{c}N_{i}^{e}\ ^{F}g_{be} & \ ^{F}g_{ab}%
\end{array}%
\right]  \label{fndm}
\end{equation}%
and satisfying the equations (for instance, for normalized flows)
\begin{eqnarray*}
\frac{\partial }{\partial \tau }\ ^{F}g_{\underline{\alpha }\underline{\beta
}} &=&-2\ _{\shortmid }^{F}R_{\underline{\alpha }\underline{\beta }}+\frac{2r%
}{5}\ ^{F}g_{\underline{\alpha }\underline{\beta }}, \\
\ ^{F}g_{\underline{\alpha }\underline{\beta }\mid _{\tau =0}} &=&\ ^{F}g_{%
\underline{\alpha }\underline{\beta }}^{[0]}(u).
\end{eqnarray*}
\end{theorem}

The Finsler--Ricci flows are distinguished from the usual (unconstrained)
flows of Riemannian metrics by existence of additional evolutions of
preferred N--adapted frames (see Proposition \ref{pddv}):

\begin{corollary}
\label{cfrf}The evolution, for all ''time'' $\tau \in \lbrack 0,\tau _{0}),$
of preferred frames on a Finsler space
\begin{equation*}
\ ^{F}\mathbf{e}_{\alpha }(\tau )=\ ^{F}\mathbf{e}_{\alpha }^{\ \underline{%
\alpha }}(\tau ,u)\partial _{\underline{\alpha }}
\end{equation*}%
is defined by the coefficients
\begin{equation}
\ ^{F}\mathbf{e}_{\alpha }^{\ \underline{\alpha }}(\tau ,u)=\left[
\begin{array}{cc}
\ ^{F}e_{i}^{\ \underline{i}}(\tau ,u) & ~^{c}N_{i}^{b}(\tau ,u)\
^{F}e_{b}^{\ \underline{a}}(\tau ,u) \\
0 & \ ^{F}e_{a}^{\ \underline{a}}(\tau ,u)%
\end{array}%
\right] ,\   \label{vt1ff}
\end{equation}%
with
\begin{equation*}
\ ^{F}g_{ij}(\tau )=\ ^{F}e_{i}^{\ \underline{i}}(\tau ,u)\ ^{F}e_{j}^{\
\underline{j}}(\tau ,u)\eta _{\underline{i}\underline{j}},
\end{equation*}%
where $\eta _{\underline{i}\underline{j}}=diag[\pm 1,...\pm 1]$ establish
the signature of $\ ^{F}g_{\underline{\alpha }\underline{\beta }}^{[0]}(u),$
is given by equations
\begin{equation}
\frac{\partial }{\partial \tau }\ ^{F}e_{\alpha }^{\ \underline{\alpha }}=\
^{F}g^{\underline{\alpha }\underline{\beta }}~_{\shortmid }^{F}R_{\underline{%
\beta }\underline{\gamma }}~\ ^{F}e_{\alpha }^{\ \underline{\gamma }},
\label{eoff}
\end{equation}%
where $\ ^{F}g^{\underline{\alpha }\underline{\beta }}$ is inverse to (\ref%
{fndm}) and $_{\shortmid }^{F}R_{\underline{\beta }\underline{\gamma }}$ is
the Ricci tensor constructed from the Levi Civita coefficients of (\ref{fndm}%
).
\end{corollary}

\begin{proof}
We have to introduce the metric and N--connection coefficients (\ref%
{sdmglsff}) and (\ref{cncff}), defined by $F(\tau ),$ into (\ref{vt1}). The
equations (\ref{eoff}) are similar to (\ref{eof}), but in our case for the
N--adapted frames (\ref{vt1ff}). $\square $
\end{proof}

\vskip3pt

We note that the evolution of the Riemann and Ricci tensors and scalar
curvature defined by the Cartan d--connection, i.e. the canonical
d--connecti\-on, $~^{F}\widehat{\mathbf{\Gamma }}_{\ \alpha \beta }^{\gamma
},$ can be extracted as in Theorem \ref{th313} when the values are redefined
for the metric (\ref{fndm}) and (\ref{fdlc}).

Finally, in this section, we conclude that the Ricci flows of Finsler
metrics can be extracted from Ricci flows of Riemannian metrics by
corresponding metric ansatz, nonholonomic constraints and deformations of
linear connections, all derived canonically from fundamental Finsler
functions.

\subsection{Ricci flows of regular Lagrange systems}

There were elaborated different approaches to geometric mechanics. Here we
note those based on formulation in terms of sympletic geometry and
generalizations \cite{marsrat,deleon1,deleon2} and in terms of generalized
Finsler, i.e. Lagrange, geomery \cite{kern,ma1,ma2}. We note that the second
approach can be also equivalently redefined as an almost Hermitian geometry
(see formulas (\ref{acs}) defining the almost complex structure) and, which
is very important for applications of the theory of anholonomic Ricci flows,
modelled as a nonholonomic Riemann manifold, see Conclusion \ref{ceqlfg}.
For regular mechanical systems, we can formulate the problem: Which
fudamental Lagrange function $L(\tau )=L(\tau ,x^{i},y^{j})$ from a class of
Lagrangians parametrized by $\tau \in \lbrack 0,\tau _{0})$ will define the
evolution of Lagrange geometry, from viewpoint of the theory of Ricci flows?
The aim of this section is to schetch the key results solving this problem.

Following the formulas from Result \ref{rap5} and the methods elaborated in
previous section \ref{ssxrl}, when $F^{2}(\tau )\rightarrow L(\tau );$ $%
^{F}h_{ij}(\tau )\rightarrow \ ^{L}g_{ij}(\tau ),$ see (\ref{ffh}) and (\ref%
{lqf}); $~^{c}N_{i}^{a}(\tau )\rightarrow $ $^{L}N_{j}^{i}(\tau ),$ see (\ref%
{cncff}) and (\ref{cncl}); $~^{c}\mathbf{g}(\tau )\rightarrow \ ^{L}\mathbf{g%
}(\tau )\mathbf{,}$ see (\ref{sdmglsff}) and (\ref{slm}); $~^{F}\widehat{%
\mathbf{\Gamma }}_{\ \beta \gamma }^{\alpha }(\tau )\rightarrow $ $\ ^{L}%
\widehat{\mathbf{\Gamma }}_{\ \beta \gamma }^{\alpha }(\tau ),$ see \ (\ref%
{fdlc}) and (\ref{cdctb}), where all values labelel by up--left $"L"$ are
canonically defined by $L(\tau ),$ we prove (generalizations of Theorem \ref%
{tfrf} and Corollary \ref{cfrf}):

\begin{theorem}
The Lagrange--Ricci flows for regular Lagrangians $L(\tau )$ can be
extracted from usual Ricci flows of Riemannian metrics paramet\-riz\-ed in
the form%
\begin{equation*}
\ ^{L}g_{\underline{\alpha }\underline{\beta }}(\tau )=\left[
\begin{array}{cc}
\ ^{L}g_{ij}+~^{L}N_{i}^{a}~^{L}N_{j}^{b}\ ^{L}g_{ab} & ~^{L}N_{j}^{e}\
^{L}g_{ae} \\
~^{L}N_{i}^{e}\ ^{L}g_{be} & \ ^{L}g_{ab}%
\end{array}%
\right]
\end{equation*}%
and satisfying the equations (for instance, normalized)
\begin{eqnarray*}
\frac{\partial }{\partial \tau }\ ^{L}g_{\underline{\alpha }\underline{\beta
}} &=&-2\ _{\shortmid }^{L}R_{\underline{\alpha }\underline{\beta }}+\frac{2r%
}{5}\ ^{L}g_{\underline{\alpha }\underline{\beta }}, \\
\ ^{L}g_{\underline{\alpha }\underline{\beta }\mid _{\tau =0}} &=&\ ^{L}g_{%
\underline{\alpha }\underline{\beta }}^{[0]}(u),
\end{eqnarray*}%
where $\ _{\shortmid }^{L}R_{\underline{\alpha }\underline{\beta }}(\tau )$
are the Ricci tensors constructed from the Levi Civita connections of
metrics $\ ^{L}g_{\underline{\alpha }\underline{\beta }}(\tau ).$
\end{theorem}

The Lagrange--Ricci flows are are characterized by the evolutions of
preferred N--adapted frames (see Proposition \ref{pddv}):

\begin{corollary}
The evolution, for all time $\tau \in \lbrack 0,\tau _{0}),$ of preferred
frames on a Lagrange space
\begin{equation*}
\ ^{L}\mathbf{e}_{\alpha }(\tau )=\ ^{L}\mathbf{e}_{\alpha }^{\ \underline{%
\alpha }}(\tau ,u)\partial _{\underline{\alpha }}
\end{equation*}%
is defined by the coefficients
\begin{equation*}
\ ^{L}\mathbf{e}_{\alpha }^{\ \underline{\alpha }}(\tau ,u)=\left[
\begin{array}{cc}
\ ^{L}e_{i}^{\ \underline{i}}(\tau ,u) & ~^{L}N_{i}^{b}(\tau ,u)\
^{L}e_{b}^{\ \underline{a}}(\tau ,u) \\
0 & \ ^{L}e_{a}^{\ \underline{a}}(\tau ,u)%
\end{array}%
\right] ,\
\end{equation*}%
with
\begin{equation*}
\ ^{L}g_{ij}(\tau )=\ ^{L}e_{i}^{\ \underline{i}}(\tau ,u)\ ^{L}e_{j}^{\
\underline{j}}(\tau ,u)\eta _{\underline{i}\underline{j}},
\end{equation*}%
where $\eta _{\underline{i}\underline{j}}=diag[\pm 1,...\pm 1]$ establish
the signature of $\ ^{L}g_{\underline{\alpha }\underline{\beta }}^{[0]}(u),$
is given by equations
\begin{equation*}
\frac{\partial }{\partial \tau }\ ^{L}e_{\alpha }^{\ \underline{\alpha }}=\
^{L}g^{\underline{\alpha }\underline{\beta }}~_{\shortmid }^{L}R_{\underline{%
\beta }\underline{\gamma }}~\ ^{L}e_{\alpha }^{\ \underline{\gamma }}.
\end{equation*}
\end{corollary}

We conclude that the Ricci flows of Lagrange metrics can be extracted from
Ricci flows of Riemannian metrics by corresponding metric ansatz,
nonholonomic constraints and deformations of linear connections, all derived
canonically for regular Lagrange functions.

\subsection{Generalized Lagrange--Ricci flows}

\label{ssxrf}We have the result that any mechanical system with a regular
Lagrangian $L(x,y)$ can be geometrized canonically in terms of nonholonomic
Riemann geometry, see Conclusion \ref{ceqlfg}, and for certain conditions
such configurations generate exact solutions of the gravitational field
equations in the Einstein gravity and/or its string/gauge generalizations,
see Result \ref{result2} and Theorem \ref{tlr}. In other turn, for any
symmetric tensor $g_{ij}=\delta _{i}^{a}$ $\delta _{j}^{b}h_{ab}(x,y)$ on a
manifold $\mathbf{V}^{n+n}$ we can generate a Lagrange space model, see
section \ref{ssgls}. The aim of this section is to show how we can construct
nonholonomic Ricci flows with effective Lagrangians starting from an
arbitrary family $g_{ij}(\tau )=\delta _{i}^{a}$ $\delta _{j}^{b}h_{ab}(\tau
,x,y).$\footnote{%
for some special cases, we can consider that $g_{ij}(\tau )$ is defined by
certain families of exact (non) holonomic solutions of the Einstein
equations or of the Ricci flow equations modelling Ricci flows of some
effective Lagrangians}

The values $h_{ab}(\tau )$ of constant signature defines a family of
absolute energies $\varepsilon (\tau )=h_{ab}(\tau ,x,y)\ y^{a}y^{b}$and
d--metrics of type (\ref{sdmgls}),%
\begin{eqnarray}
~^{\varepsilon }\mathbf{g}(\tau ) &=&\ h_{ij}(\tau ,x,y)\ \left(
e^{i}\otimes e^{j}+\ ~^{\varepsilon }\mathbf{e}^{i}(\tau )\otimes
~^{\varepsilon }\mathbf{e}^{j}(\tau )\right) ,  \notag \\
~^{\varepsilon }\mathbf{e}^{i}(\tau ) &=&dy^{i}+~^{\varepsilon
}N_{i}^{a}(\tau ,x,y)dx^{i},  \label{dmrfel}
\end{eqnarray}%
where the $\tau $--parametrized N--connection coefficients
\begin{equation}
~^{\varepsilon }N_{i}^{a}(\tau ,x,y)=\frac{\partial ~^{\varepsilon
}G^{a}(\tau )}{\partial y^{i}},  \label{cncelrf}
\end{equation}%
with%
\begin{equation*}
~^{\varepsilon }G^{a}(\tau )=\frac{1}{2}~^{\varepsilon }\widetilde{h}%
^{ab}(\tau )\left( y^{k}\frac{\partial ^{2}\varepsilon (\tau )}{\partial
y^{b}\partial x^{k}}-\delta _{b}^{k}\frac{\partial \varepsilon (\tau )}{%
\partial x^{k}}\right) ,
\end{equation*}%
are defined for nondegenerated Hessians
\begin{equation}
~^{\varepsilon }\widetilde{h}_{ab}(\tau )=\frac{1}{2}\frac{\partial
^{2}\varepsilon (\tau )}{\partial y^{a}\partial y^{b}},  \label{hesselrf}
\end{equation}%
when $\det |\widetilde{h}|\neq 0.$

For any fixed value of $\tau ,$ the existence of fundamental geometric
objects (\ref{dmrfel}), (\ref{cncelrf}) and (\ref{hesselrf}) follows from
Theorem \ref{tcsl1}. Similarly, the Theorem \ref{tcslg} states a modelling
by $h_{ab}(\tau )$ of families of Lagrange spaces enabled with canonical
N--connections $\ ^{\varepsilon }\mathbf{N(\tau ),}$ almost complex
structure $^{c}\mathbf{F(\tau ),}$ d--metrics $~^{c}\mathbf{g}(\tau )$ and
d--connections $~^{c}\widehat{\mathbf{D}}(\tau )$ structures defined
respectively by effective regular Lagrangians $~^{\varepsilon }L(\tau ,x,y)=%
\sqrt{|\varepsilon (\tau ,x,y)|}$ and theirs Hessians $~^{\varepsilon }%
\widetilde{h}_{ab}(\tau ,x,y)$ (\ref{hesselrf}). The results of previous
section \ref{ssxrf} can be reformulated in the form (with proofs being
similar for those for Theorem \ref{tfrf} and Corollary \ref{cfrf}, but with $%
~^{\varepsilon }L$ instead of $F^{2}$ and $~^{\varepsilon }N_{i}^{a}$
instead of $~^{c}N_{i}^{a},...$):

\begin{theorem}
The generalized Lagrange--Ricci flows for regular effective Lagrangians $%
~^{\varepsilon }L(\tau )$ $\ $derived from a family of symmetric tensors $%
h_{ab}(\tau ,x,y)$ can be extracted from usual Ricci flows of Riemannian
metrics paramet\-riz\-ed in the form%
\begin{equation*}
\ ^{\varepsilon }g_{\underline{\alpha }\underline{\beta }}(\tau )=\left[
\begin{array}{cc}
~^{\varepsilon }\widetilde{h}_{ij}+~^{\varepsilon }N_{i}^{a}~^{\varepsilon
}N_{j}^{b}\ ~^{\varepsilon }\widetilde{h}_{ab} & ~^{\varepsilon
}N_{j}^{e}~^{\varepsilon }\widetilde{h}_{ae} \\
~^{\varepsilon }N_{i}^{e}~^{\varepsilon }\widetilde{h}_{be} & ~^{\varepsilon
}\widetilde{h}_{ab}%
\end{array}%
\right]
\end{equation*}%
and satisfying the equations (for instance, normalized)
\begin{eqnarray*}
\frac{\partial }{\partial \tau }\ ^{\varepsilon }g_{\underline{\alpha }%
\underline{\beta }} &=&-2\ _{\shortmid }^{\varepsilon }R_{\underline{\alpha }%
\underline{\beta }}+\frac{2r}{5}\ ^{\varepsilon }g_{\underline{\alpha }%
\underline{\beta }}, \\
\ ^{\varepsilon }g_{\underline{\alpha }\underline{\beta }\mid _{\tau =0}}
&=&\ ^{\varepsilon }g_{\underline{\alpha }\underline{\beta }}^{[0]}(u),
\end{eqnarray*}%
where $\ _{\shortmid }^{\varepsilon }R_{\underline{\alpha }\underline{\beta }%
}(\tau )$ are the Ricci tensors constructed from the Levi Civita connections
of metrics $\ ^{\varepsilon }g_{\underline{\alpha }\underline{\beta }}(\tau
).$
\end{theorem}

The evolutions of preferred N--adapted frames (see Proposition \ref{pddv})
defined by generalized Lagrange--Ricci flows is stated by

\begin{corollary}
The evolution, for all time $\tau \in \lbrack 0,\tau _{0}),$ of preferred
frames on an effective Lagrange space
\begin{equation*}
~^{\varepsilon }\mathbf{e}_{\alpha }(\tau )=~^{\varepsilon }\mathbf{e}%
_{\alpha }^{\ \underline{\alpha }}(\tau ,u)\partial _{\underline{\alpha }}
\end{equation*}%
is defined by the coefficients
\begin{equation*}
\ ~^{\varepsilon }\mathbf{e}_{\alpha }^{\ \underline{\alpha }}(\tau ,u)=%
\left[
\begin{array}{cc}
~^{\varepsilon }e_{i}^{\ \underline{i}}(\tau ,u) & ~^{\varepsilon
}N_{i}^{b}(\tau ,u)~^{\varepsilon }e_{b}^{\ \underline{a}}(\tau ,u) \\
0 & ~^{\varepsilon }e_{a}^{\ \underline{a}}(\tau ,u)%
\end{array}%
\right] ,\
\end{equation*}%
with
\begin{equation*}
~^{\varepsilon }\widetilde{h}_{ij}(\tau )=~^{\varepsilon }e_{i}^{\
\underline{i}}(\tau ,u)~^{\varepsilon }e_{j}^{\ \underline{j}}(\tau ,u)\eta
_{\underline{i}\underline{j}},
\end{equation*}%
where $\eta _{\underline{i}\underline{j}}=diag[\pm 1,...\pm 1]$ establish
the signature of $\ ~^{\varepsilon }g_{\underline{\alpha }\underline{\beta }%
}^{[0]}(u),$ is given by equations
\begin{equation*}
\frac{\partial }{\partial \tau }~^{\varepsilon }e_{\alpha }^{\ \underline{%
\alpha }}=~^{\varepsilon }g^{\underline{\alpha }\underline{\beta }%
}~_{\shortmid }^{\varepsilon }R_{\underline{\beta }\underline{\gamma }}~\
~^{\varepsilon }e_{\alpha }^{\ \underline{\gamma }}.
\end{equation*}
\end{corollary}

The idea to consider absolute energies $\varepsilon $ for arbitrary
d--metrics $g_{ij}(x,y),$ in order to define effective (generalized)
Lagrange spaces, was proposed in Refs.\cite{ma1,ma2}. In Introduction and
Part I of monograph \cite{vsgg}, it was proven that certain type of
gravitational interactions can be modelled as generalized Lagrange--Finsler
geometries and inversely, certain classes of generalized Finsler geometries
can be modelled on N--anholonomic manifolds, even as exact solutions of
gravitational field equations. The approach elaborated by Romanian geometers
and physicists \cite{ma1,ma2,mhh,vncg,vsgg,vcla} originates from G.
Vranceanu and Z. Horac works \cite{vr1,hor} on nonholonomic manifolds and
mechanical systems, see a review of results and recent developments in Ref. %
\cite{bejf}. Recently, there were proposed various models of ''analogous
gravity'', see a review in Ref. \cite{blv}, which do not apply the methods
of Finsler geometry and the formalism of nonlinear connections.

\vskip5pt

\textbf{Acknowledgement: } The work is performed during a visit at the
Fields Institute. The author thanks M. Anastasiei, D. Hrimiuc, D. Singleton
and M. Visinescu for substantial support. He is also grateful to V. Ob\v{a}%
deanu and A. Bejancu for very important references on the geometry of
nonholonomic manifolds.

\setcounter{equation}{0} \renewcommand{\theequation}
{A.\arabic{equation}} \setcounter{subsection}{0}
\renewcommand{\thesubsection}
{A.\arabic{subsection}}

\appendix

\section{Local Geometry of N--anholonomic Manifolds}

Let us consider metric structure on N--anholonomic manifold $\mathbf{V},$%
\begin{equation}
\ \breve{g}=\underline{g}_{\alpha \beta }\left( u\right) du^{\alpha }\otimes
du^{\beta }  \label{metr}
\end{equation}%
defined with respect to a local coordinate basis $du^{\alpha }=\left(
dx^{i},dy^{a}\right) $ by coefficients%
\begin{equation}
\underline{g}_{\alpha \beta }=\left[
\begin{array}{cc}
g_{ij}+N_{i}^{a}N_{j}^{b}h_{ab} & N_{j}^{e}h_{ae} \\
N_{i}^{e}h_{be} & h_{ab}%
\end{array}%
\right] .  \label{ansatz}
\end{equation}%
Such a metric (\ref{ansatz})\ is generic off--diagonal, i.e. it can not be
diagonalized by coordinate transforms if $N_{i}^{a}(u)$ are any general
functions. The condition (\ref{algn01}), for $hX\rightarrow e_{i}$ and $\
vY\rightarrow e_{a},$ transform into
\begin{equation}
\breve{g}(e_{i},e_{a})=0,\mbox{ equivalently }\underline{g}%
_{ia}-N_{i}^{b}h_{ab}=0,  \label{aux1a}
\end{equation}%
where $\underline{g}_{ia}$ $\doteqdot g(\partial /\partial x^{i},\partial
/\partial y^{a}),$ which allows us to define in a unique form the
coefficients $N_{i}^{b}=h^{ab}\underline{g}_{ia}$ where $h^{ab}$ is inverse
to $h_{ab}.$ We can write the metric $\breve{g}$ with ansatz (\ref{ansatz})\
in equivalent form, as a d--metric (\ref{m1}) adapted to a N--connection
structure, see Definition \ref{ddm}, if we define $g_{ij}\doteqdot \mathbf{g}%
\left( e_{i},e_{j}\right) $ and $h_{ab}\doteqdot \mathbf{g}\left(
e_{a},e_{b}\right) $ \ and consider the vielbeins $\mathbf{e}_{\alpha }$ and
$\mathbf{e}^{\alpha }$ to be respectively of type (\ref{dder}) and (\ref%
{ddif}).

We can say that the metric $\ \breve{g}$ (\ref{metr}) is equivalently
transformed into (\ref{m1}) by performing a frame (vielbein) transform
\begin{equation*}
\mathbf{e}_{\alpha }=\mathbf{e}_{\alpha }^{\ \underline{\alpha }}\partial _{%
\underline{\alpha }}\mbox{ and }\mathbf{e}_{\ }^{\beta }=\mathbf{e}_{\
\underline{\beta }}^{\beta }du^{\underline{\beta }}.
\end{equation*}%
with coefficients

\begin{eqnarray}
\mathbf{e}_{\alpha }^{\ \underline{\alpha }}(u) &=&\left[
\begin{array}{cc}
e_{i}^{\ \underline{i}}(u) & N_{i}^{b}(u)e_{b}^{\ \underline{a}}(u) \\
0 & e_{a}^{\ \underline{a}}(u)%
\end{array}%
\right] ,  \label{vt1} \\
\mathbf{e}_{\ \underline{\beta }}^{\beta }(u) &=&\left[
\begin{array}{cc}
e_{\ \underline{i}}^{i\ }(u) & -N_{k}^{b}(u)e_{\ \underline{i}}^{k\ }(u) \\
0 & e_{\ \underline{a}}^{a\ }(u)%
\end{array}%
\right] ,  \label{vt2}
\end{eqnarray}%
being linear on $N_{i}^{a}.$ We can consider that a N--anholonomic manifold $%
\mathbf{V}$ provided with metric structure $\breve{g}$ (\ref{metr})
(equivalently, with d--metric (\ref{m1})) is a special type of a manifold
provided with a global splitting into conventional ``horizontal'' and
``vertical'' subspaces (\ref{whitney}) induced by the ``off--diagonal''
terms $N_{i}^{b}(u)$ and a prescribed type of nonholonomic frame structure (%
\ref{anhrel}).

The N--adapted components $\mathbf{\Gamma }_{\ \beta \gamma }^{\alpha }$ of
a d--connection $\mathbf{D}_{\alpha }=(\mathbf{e}_{\alpha }\rfloor \mathbf{D}%
),$ where ''$\rfloor $'' denotes the interior product, are defined by the
equations
\begin{equation}
\mathbf{D}_{\alpha }\mathbf{e}_{\beta }=\mathbf{\Gamma }_{\ \alpha \beta
}^{\gamma }\mathbf{e}_{\gamma },\mbox{\ or \ }\mathbf{\Gamma }_{\ \alpha
\beta }^{\gamma }\left( u\right) =\left( \mathbf{D}_{\alpha }\mathbf{e}%
_{\beta }\right) \rfloor \mathbf{e}^{\gamma }.  \label{dcon1}
\end{equation}%
The N--adapted splitting into h-- and v--covariant derivatives is stated by
\begin{equation*}
h\mathbf{D}=\{\mathbf{D}_{k}=\left( L_{jk}^{i},L_{bk\;}^{a}\right) \},%
\mbox{
and }\ v\mathbf{D}=\{\mathbf{D}_{c}=\left( C_{jc}^{i},C_{bc}^{a}\right) \},
\end{equation*}%
where, by definition,
\begin{equation*}
L_{jk}^{i}=\left( \mathbf{D}_{k}\mathbf{e}_{j}\right) \rfloor e^{i},\quad
L_{bk}^{a}=\left( \mathbf{D}_{k}e_{b}\right) \rfloor \mathbf{e}%
^{a},~C_{jc}^{i}=\left( \mathbf{D}_{c}\mathbf{e}_{j}\right) \rfloor
e^{i},\quad C_{bc}^{a}=\left( \mathbf{D}_{c}e_{b}\right) \rfloor \mathbf{e}%
^{a}.
\end{equation*}%
The components $\mathbf{\Gamma }_{\ \alpha \beta }^{\gamma }=\left(
L_{jk}^{i},L_{bk}^{a},C_{jc}^{i},C_{bc}^{a}\right) $ completely define a
d--connection $\mathbf{D}$ on a N--anholonomic manifold $\mathbf{V}.$

The simplest way to perform computations with d--connections is to use
N--adapted differential forms like
\begin{equation}
\mathbf{\Gamma }_{\ \beta }^{\alpha }=\mathbf{\Gamma }_{\ \beta \gamma
}^{\alpha }\mathbf{e}^{\gamma }  \label{dconf}
\end{equation}%
with the coefficients defined with respect to (\ref{ddif}) and (\ref{dder}).
For instance, torsion be computed in the form
\begin{equation}
\mathcal{T}^{\alpha }\doteqdot \mathbf{De}^{\alpha }=d\mathbf{e}^{\alpha
}+\Gamma _{\ \beta }^{\alpha }\wedge \mathbf{e}^{\beta }.  \label{tors}
\end{equation}%
Locally it is characterized by (N--adapted) d--torsion coefficients
\begin{eqnarray}
T_{\ jk}^{i} &=&L_{\ jk}^{i}-L_{\ kj}^{i},\ T_{\ ja}^{i}=-T_{\ aj}^{i}=C_{\
ja}^{i},\ T_{\ ji}^{a}=\Omega _{\ ji}^{a},\   \notag \\
T_{\ bi}^{a} &=&-T_{\ ib}^{a}=\frac{\partial N_{i}^{a}}{\partial y^{b}}-L_{\
bi}^{a},\ T_{\ bc}^{a}=C_{\ bc}^{a}-C_{\ cb}^{a}.  \label{dtors}
\end{eqnarray}

By a straightforward d--form calculus, we can find the N--adapted components
of the curvature
\begin{equation}
\mathcal{R}_{~\beta }^{\alpha }\doteqdot \mathbf{D\Gamma }_{\ \beta
}^{\alpha }=d\mathbf{\Gamma }_{\ \beta }^{\alpha }-\mathbf{\Gamma }_{\ \beta
}^{\gamma }\wedge \mathbf{\Gamma }_{\ \gamma }^{\alpha }=\mathbf{R}_{\ \beta
\gamma \delta }^{\alpha }\mathbf{e}^{\gamma }\wedge \mathbf{e}^{\delta },
\label{curv}
\end{equation}%
of a d--connection $\mathbf{D},$ i.e. the d--curvatures from Theorem \ref%
{thr}:
\begin{eqnarray}
R_{\ hjk}^{i} &=&e_{k}L_{\ hj}^{i}-e_{j}L_{\ hk}^{i}+L_{\ hj}^{m}L_{\
mk}^{i}-L_{\ hk}^{m}L_{\ mj}^{i}-C_{\ ha}^{i}\Omega _{\ kj}^{a},  \notag \\
R_{\ bjk}^{a} &=&e_{k}L_{\ bj}^{a}-e_{j}L_{\ bk}^{a}+L_{\ bj}^{c}L_{\
ck}^{a}-L_{\ bk}^{c}L_{\ cj}^{a}-C_{\ bc}^{a}\Omega _{\ kj}^{c},  \notag \\
R_{\ jka}^{i} &=&e_{a}L_{\ jk}^{i}-D_{k}C_{\ ja}^{i}+C_{\ jb}^{i}T_{\
ka}^{b},  \label{dcurv} \\
R_{\ bka}^{c} &=&e_{a}L_{\ bk}^{c}-D_{k}C_{\ ba}^{c}+C_{\ bd}^{c}T_{\
ka}^{c},  \notag \\
R_{\ jbc}^{i} &=&e_{c}C_{\ jb}^{i}-e_{b}C_{\ jc}^{i}+C_{\ jb}^{h}C_{\
hc}^{i}-C_{\ jc}^{h}C_{\ hb}^{i},  \notag \\
R_{\ bcd}^{a} &=&e_{d}C_{\ bc}^{a}-e_{c}C_{\ bd}^{a}+C_{\ bc}^{e}C_{\
ed}^{a}-C_{\ bd}^{e}C_{\ ec}^{a}.  \notag
\end{eqnarray}

Contracting respectively the components of (\ref{dcurv}), one proves that
the Ricci tensor $\mathbf{R}_{\alpha \beta }\doteqdot \mathbf{R}_{\ \alpha
\beta \tau }^{\tau }$ is characterized by h- v--components, i.e. d--tensors,%
\begin{equation}
R_{ij}\doteqdot R_{\ ijk}^{k},\ \ R_{ia}\doteqdot -R_{\ ika}^{k},\
R_{ai}\doteqdot R_{\ aib}^{b},\ R_{ab}\doteqdot R_{\ abc}^{c}.
\label{dricci}
\end{equation}%
It should be noted that this tensor is not symmetric for arbitrary
d--connecti\-ons $\mathbf{D}.$

The scalar curvature of a d--connection is
\begin{equation}
\ ^{s}\mathbf{R}\doteqdot \mathbf{g}^{\alpha \beta }\mathbf{R}_{\alpha \beta
}=g^{ij}R_{ij}+h^{ab}R_{ab},  \label{sdccurv}
\end{equation}%
defined by a sum the h-- and v--components of (\ref{dricci}) and d--metric (%
\ref{m1}).

The Einstein tensor is defined and computed in standard form
\begin{equation}
\mathbf{G}_{\alpha \beta }=\mathbf{R}_{\alpha \beta }-\frac{1}{2}\mathbf{g}%
_{\alpha \beta }\ ^{s}\mathbf{R}  \label{enstdt}
\end{equation}

One exists a minimal extension of the Levi Civita connection $\nabla $ to a
canonical d--connection $\widehat{\mathbf{D}}$ which is defined only by a
metric $\breve{g}$ is metric compatible, with $\widehat{T}_{\ jk}^{i}=0$ and
$\widehat{T}_{\ bc}^{a}=0$ but $\widehat{T}_{\ ja}^{i},\widehat{T}_{\
ji}^{a} $ and $\widehat{T}_{\ bi}^{a}$ are not zero, see (\ref{dtors}). The
coefficient $\widehat{\mathbf{\Gamma }}_{\ \alpha \beta }^{\gamma }=\left(
\widehat{L}_{jk}^{i},\widehat{L}_{bk}^{a},\widehat{C}_{jc}^{i},\widehat{C}%
_{bc}^{a}\right) $ of this connection, with respect to the N--adapted
frames, are defined are computed:
\begin{eqnarray}
\widehat{L}_{jk}^{i} &=&\frac{1}{2}g^{ir}\left(
e_{k}g_{jr}+e_{j}g_{kr}-e_{r}g_{jk}\right) ,  \label{candcon} \\
\widehat{L}_{bk}^{a} &=&e_{b}(N_{k}^{a})+\frac{1}{2}h^{ac}\left(
e_{k}h_{bc}-h_{dc}\ e_{b}N_{k}^{d}-h_{db}\ e_{c}N_{k}^{d}\right) ,  \notag \\
\widehat{C}_{jc}^{i} &=&\frac{1}{2}g^{ik}e_{c}g_{jk},\ \widehat{C}_{bc}^{a}=%
\frac{1}{2}h^{ad}\left( e_{c}h_{bd}+e_{c}h_{cd}-e_{d}h_{bc}\right) .  \notag
\end{eqnarray}

The Levi Civita linear connection $\bigtriangledown =\{\ _{\shortmid }\Gamma
_{\beta \gamma }^{\alpha }\},$ uniquely defined by the conditions $~\
_{\shortmid }\mathcal{T}=0$ and $\bigtriangledown g=0,$ is not adapted to
the distribution (\ref{whitney}). Let us parametrize the coefficients in the
form
\begin{equation*}
_{\shortmid }\Gamma _{\beta \gamma }^{\alpha }=\left( _{\shortmid
}L_{jk}^{i},_{\shortmid }L_{jk}^{a},_{\shortmid }L_{bk}^{i},\ _{\shortmid
}L_{bk}^{a},_{\shortmid }C_{jb}^{i},_{\shortmid }C_{jb}^{a},_{\shortmid
}C_{bc}^{i},_{\shortmid }C_{bc}^{a}\right) ,
\end{equation*}%
where
\begin{eqnarray*}
\bigtriangledown _{\mathbf{e}_{k}}(\mathbf{e}_{j}) &=&_{\shortmid }L_{jk}^{i}%
\mathbf{e}_{i}+_{\shortmid }L_{jk}^{a}e_{a},\bigtriangledown _{\mathbf{e}%
_{k}}(e_{b})=_{\shortmid }L_{bk}^{i}\mathbf{e}_{i}+\ _{\shortmid
}L_{bk}^{a}e_{a}, \\
\bigtriangledown _{e_{b}}(\mathbf{e}_{j}) &=&_{\shortmid }C_{jb}^{i}\mathbf{e%
}_{i}+_{\shortmid }C_{jb}^{a}e_{a},\bigtriangledown
_{e_{c}}(e_{b})=_{\shortmid }C_{bc}^{i}\mathbf{e}_{i}+_{\shortmid
}C_{bc}^{a}e_{a}.
\end{eqnarray*}%
A straightforward calculus\footnote{%
Such results were originally considered by R. Miron and M. Anastasiei for
vector bundles provided with N--connection and metric structures, see Ref. %
\cite{ma2}. Similar proofs hold true for any nonholonomic manifold provided
with a prescribed N--connection structures.} shows that the coefficients of
the Levi-Civita connection can be expressed in the form%
\begin{eqnarray}
\ _{\shortmid }L_{jk}^{i} &=&L_{jk}^{i},\ _{\shortmid
}L_{jk}^{a}=-C_{jb}^{i}g_{ik}h^{ab}-\frac{1}{2}\Omega _{jk}^{a},
\label{lccon} \\
\ _{\shortmid }L_{bk}^{i} &=&\frac{1}{2}\Omega _{jk}^{c}h_{cb}g^{ji}-\frac{1%
}{2}(\delta _{j}^{i}\delta _{k}^{h}-g_{jk}g^{ih})C_{hb}^{j},  \notag \\
\ _{\shortmid }L_{bk}^{a} &=&L_{bk}^{a}+\frac{1}{2}(\delta _{c}^{a}\delta
_{d}^{b}+h_{cd}h^{ab})\left[ L_{bk}^{c}-e_{b}(N_{k}^{c})\right] ,  \notag \\
\ _{\shortmid }C_{kb}^{i} &=&C_{kb}^{i}+\frac{1}{2}\Omega
_{jk}^{a}h_{cb}g^{ji}+\frac{1}{2}(\delta _{j}^{i}\delta
_{k}^{h}-g_{jk}g^{ih})C_{hb}^{j},  \notag \\
\ _{\shortmid }C_{jb}^{a} &=&-\frac{1}{2}(\delta _{c}^{a}\delta
_{b}^{d}-h_{cb}h^{ad})\left[ L_{dj}^{c}-e_{d}(N_{j}^{c})\right] ,\
_{\shortmid }C_{bc}^{a}=C_{bc}^{a},  \notag \\
\ _{\shortmid }C_{ab}^{i} &=&-\frac{g^{ij}}{2}\left\{ \left[
L_{aj}^{c}-e_{a}(N_{j}^{c})\right] h_{cb}+\left[ L_{bj}^{c}-e_{b}(N_{j}^{c})%
\right] h_{ca}\right\} ,  \notag
\end{eqnarray}%
where $\Omega _{jk}^{a}$ are computed as in formula (\ref{ncurv}). For
certain considerations, it is convenient to express
\begin{equation}
\ _{\shortmid }\Gamma _{\ \alpha \beta }^{\gamma }=\widehat{\mathbf{\Gamma }}%
_{\ \alpha \beta }^{\gamma }+\ _{\shortmid }Z_{\ \alpha \beta }^{\gamma }
\label{cdeft}
\end{equation}%
where the explicit components of distorsion tensor $\ _{\shortmid }Z_{\
\alpha \beta }^{\gamma }$ can be defined by comparing the formulas (\ref%
{lccon}) and (\ref{candcon}):%
\begin{eqnarray}
\ _{\shortmid }Z_{jk}^{i} &=&0,\ _{\shortmid
}Z_{jk}^{a}=-C_{jb}^{i}g_{ik}h^{ab}-\frac{1}{2}\Omega _{jk}^{a},  \notag \\
\ _{\shortmid }Z_{bk}^{i} &=&\frac{1}{2}\Omega _{jk}^{c}h_{cb}g^{ji}-\frac{1%
}{2}(\delta _{j}^{i}\delta _{k}^{h}-g_{jk}g^{ih})C_{hb}^{j},  \notag \\
\ _{\shortmid }Z_{bk}^{a} &=&\frac{1}{2}(\delta _{c}^{a}\delta
_{d}^{b}+h_{cd}h^{ab})\left[ L_{bk}^{c}-e_{b}(N_{k}^{c})\right] ,  \notag \\
\ _{\shortmid }Z_{kb}^{i} &=&\frac{1}{2}\Omega _{jk}^{a}h_{cb}g^{ji}+\frac{1%
}{2}(\delta _{j}^{i}\delta _{k}^{h}-g_{jk}g^{ih})C_{hb}^{j},  \notag \\
\ _{\shortmid }Z_{jb}^{a} &=&-\frac{1}{2}(\delta _{c}^{a}\delta
_{b}^{d}-h_{cb}h^{ad})\left[ L_{dj}^{c}-e_{d}(N_{j}^{c})\right] ,\
_{\shortmid }Z_{bc}^{a}=0,  \label{cdeftc} \\
\ _{\shortmid }Z_{ab}^{i} &=&-\frac{g^{ij}}{2}\left\{ \left[
L_{aj}^{c}-e_{a}(N_{j}^{c})\right] h_{cb}+\left[ L_{bj}^{c}-e_{b}(N_{j}^{c})%
\right] h_{ca}\right\} .  \notag
\end{eqnarray}%
It should be emphasized that all components of $\ _{\shortmid }\Gamma _{\
\alpha \beta }^{\gamma },\widehat{\mathbf{\Gamma }}_{\ \alpha \beta
}^{\gamma }$ and$\ _{\shortmid }Z_{\ \alpha \beta }^{\gamma }$ are defined
by the coefficients of d--metric (\ref{m1}) and N--connection (\ref{coeffnc}%
), or equivalently by the coefficients of the corresponding generic
off--diagonal metric\ (\ref{ansatz}).

For a differentiable Lagrangian $L(x,y),$ i.e. a fundamental Lagrange
function, is defined by a map $L:(x,y)\in TM\rightarrow L(x,y)\in \mathbb{R}$
of class $\mathcal{C}^{\infty }$ on $\widetilde{TM}=TM\backprime \{0\}$ and
continuous on the null section $0:\ M\rightarrow TM$ of $\pi $ one have been %
\cite{ma1,ma2} established the following results:

\begin{result}
\begin{enumerate}
\item \label{rap5}The Euler--Lagrange equations%
\begin{equation*}
\frac{d}{d\tau }\left( \frac{\partial L}{\partial y^{i}}\right) -\frac{%
\partial L}{\partial x^{i}}=0
\end{equation*}%
where $y^{i}=\frac{dx^{i}}{d\varsigma }$ for $x^{i}(\varsigma )$ depending
on parameter $\varsigma ,$ are equivalent to the ``nonlinear'' geodesic
equations
\begin{equation*}
\frac{d^{2}x^{i}}{d\tau ^{2}}+2G^{i}(x^{k},\frac{dx^{j}}{d\varsigma })=0
\end{equation*}%
defining paths of a canonical semispray%
\begin{equation*}
S=y^{i}\frac{\partial }{\partial x^{i}}-2G^{i}(x,y)\frac{\partial }{\partial
y^{i}}
\end{equation*}%
where
\begin{equation*}
2G^{i}(x,y)=\frac{1}{2}\ ^{L}g^{ij}\left( \frac{\partial ^{2}L}{\partial
y^{i}\partial x^{k}}y^{k}-\frac{\partial L}{\partial x^{i}}\right)
\end{equation*}%
with $^{L}g^{ij}$ being inverse to (\ref{lqf}).

\item There exists on $\mathbf{V\simeq }$ $\widetilde{TM}$ a canonical
N--connection $\ $%
\begin{equation}
^{L}N_{j}^{i}=\frac{\partial G^{i}(x,y)}{\partial y^{i}}  \label{cncl}
\end{equation}%
defined by the fundamental Lagrange function $L(x,y),$ which prescribes
nonholonomic frame structures of type (\ref{dder}) and (\ref{ddif}), $^{L}%
\mathbf{e}_{\nu }=(\mathbf{e}_{i},e_{a})$ and $^{L}\mathbf{e}^{\mu }=(e^{i},%
\mathbf{e}^{a}).$

\item There is a canonical metric structure%
\begin{equation}
\ ^{L}\mathbf{g}=\ g_{ij}(x,y)\ e^{i}\otimes e^{j}+\ g_{ij}(x,y)\ \mathbf{e}%
^{i}\otimes \mathbf{e}^{j}  \label{slm}
\end{equation}%
constructed as a Sasaki type lift from $M$ \ for $g_{ij}(x,y)=\
^{L}g_{ij}(x,y),$ see (\ref{lqf}).

\item There is a unique metrical and, in this case, torsionless canonical
d--connection$\ ^{L}\mathbf{D}=(hD,vD)$ with the nontrivial coefficients
with respect to $^{L}\mathbf{e}_{\nu }$ and $^{L}\mathbf{e}^{\mu }$
paramet\-riz\-ed respectively $\ ^{L}\widehat{\mathbf{\Gamma }}_{\ \beta
\gamma }^{\alpha }=(\widehat{L}_{\ jk}^{i},\widehat{C}_{bc}^{a}),$ for
\begin{eqnarray}
\widehat{L}_{\ jk}^{i} &=&\frac{1}{2}g^{ih}(\mathbf{e}_{k}g_{jh}+\mathbf{e}%
_{j}g_{kh}-\mathbf{e}_{h}g_{jk}),  \label{cdctb} \\
\widehat{C}_{\ jk}^{i} &=&\frac{1}{2}%
g^{ih}(e_{k}g_{jh}+e_{j}g_{kh}-e_{h}h_{jk})  \notag
\end{eqnarray}%
defining the generalized Christoffel symbols, where (for simplicity, we
omitted the left up labels $(L)$ for N--adapted bases).
\end{enumerate}
\end{result}

We conclude that any regular Lagrange mechanics can be geometrized as a
nonholonomic Riemann manifold $\mathbf{V}$ equipped with canonical
N--connection (\ref{cncl}) and adapted d--connection (\ref{cdctb}) and
d--metric structures (\ref{slm}) all induced by a $L(x,y).$

Now we show how N--anholonomic configurations can defined in gravity
theories. In this case, it is convenient to work on a general manifold $%
\mathbf{V},\dim \mathbf{V}=n+m$ enabled with a global N--connection
structure, instead of the tangent bundle $\widetilde{TM}.$

We summarize here some geometric properties of gravitational models with
nontrivial N--anholonomic structure \cite{vncg,vsgg}.

\begin{result}
\label{result2}Various classes of vacuum and nonvacuum exact solutions of (%
\ref{einsta}) parametrized by generic off--diagonal metrics, nonholonomic
vielbeins and Levi Civita or non--Riemannian connections in Einstein and
extra dimension gravity models define explicit examples of N--anholonomic
Einstein--Cartan (in particular, Einstein) spaces.
\end{result}

It should be noted that a subclass of N--anholonomic Einstein spaces was
related to generic off--diagonal solutions in general relativity by such
nonholonomic constraints when $\mathbf{Ric}(\widehat{\mathbf{D}})=Ric(%
\mathbf{\nabla })$ even $\widehat{\mathbf{D}}\neq \nabla ,$ where $\widehat{%
\mathbf{D}}$ is the canonical d--connection and $\nabla $ is the
Levi--Civita connection.

A direction in modern gravity is connected to analogous gravity models when
certain gravitational effects and, for instance, black hole configurations
are modelled by optical and acoustic media, see a recent review or results
in \cite{blv}. Following our approach on geometric unification of gravity
and Lagrange regular mechanics in terms of N--anholonomic spaces, one holds

\begin{theorem}
\label{tlr}A Lagrange (Finsler) space can be canonically modelled as an
exact solution of the Einstein equations (\ref{einsta}) on a N--anholonomic
Riemann--Cartan space if and only if the canonical N--connection $\ ^{L}%
\mathbf{N} $ ($\ ^{F}\mathbf{N}$)$\mathbf{,}$ d--metric $\ ^{L}\mathbf{g}$ ($%
\ ^{F}\mathbf{g)}$ and d--connection $\ ^{L}\widehat{\mathbf{D}}$ ($\ ^{F}%
\widehat{\mathbf{D}})$ $\ $structures defined by the corresponding
fundamental Lagrange function $L(\mathbf{x,y})$ (Finsler function $F(\mathbf{%
x,y}))$ satisfy the gravitational field equations for certain physically
reasonable sources $\mathbf{\Upsilon .}$
\end{theorem}

\begin{proof}
It can be performed in local form by considering the Einstein tensor (\ref%
{enstdt}) defined by the $\ ^{L}\mathbf{N}$ ($\ ^{F}\mathbf{N}$) in the form
(\ref{cncl}) and $\ ^{L}\mathbf{g}$ ($\ ^{F}\mathbf{g)}$ in the form (\ref%
{slm}) \ inducing the canonical d--connection $\ ^{L}\widehat{\mathbf{D}}$ ($%
\ ^{F}\widehat{\mathbf{D}}).$ For certain zero or nonzero $\mathbf{\Upsilon }
$, such N--anholonomic configurations may be defined by exact solutions of
the Einstein equations for a d--connection structure. $\square $
\end{proof}

\end{document}